\documentclass[12pt,a4paper]{article}

\usepackage{color,graphics}
\usepackage[T1]{fontenc}
\usepackage{lastpage}
\usepackage{exscale}

\usepackage{lmodern}
\usepackage{bm}
\usepackage{amsthm}
\usepackage{multirow}
\usepackage{graphicx}   	
\usepackage{caption}
\usepackage{diagbox}
\usepackage[english]{babel}	
\usepackage{amsmath}
\usepackage{amssymb}
\usepackage[utf8]{inputenc}
\usepackage[dvipsnames]{xcolor}
\usepackage{subfig}
\usepackage{mathtools}
\usepackage[linesnumbered,ruled,vlined]{algorithm2e}
\usepackage{algorithmic}
\usepackage{hyperref}
\usepackage{times}
\usepackage{textcomp}
\usepackage[most]{tcolorbox}
\usepackage{physics}
\usepackage[nocompress]{cite}
\usepackage{cleveref}
\usepackage{anyfontsize}

\newlength\myindent
\newcommand\bindent{%
  \begingroup
  \setlength{\itemindent}{\myindent}
  \addtolength{\algorithmicindent}{\myindent}
}
\newcommand\eindent{\endgroup}

\newcommand{\scat}
{\mu'_{\rm s}}
\newcommand{\att}
{\mu_{\rm a}}

\newcommand{\mua}{\mu_{{\mathrm a}}}  
\newcommand{\mus}{\mu'_{{\rm s}}}

\newcommand{\revision}[1]{#1}
\newtheorem*{remark}{Remark}

\title{Sparsity promoting reconstructions via hierarchical prior models in diffuse optical tomography}
\author{\fontsize{13}{14.6}{\selectfont Anssi Manninen$^1$, Meghdoot Mozumder$^2$, Tanja Tarvainen$^{2,3}$, Andreas Hauptmann$^{1,3}$} \\ \small{$^1$Research Unit of Mathematical Sciences, University of Oulu, Oulu 90014, Finland}\\
\small{$^2$Department of Applied Physics, University of Eastern Finland, Kuopio 70211, Finland}\\
\small{$^3$Department of Computer Science, University College London, London WC1E 6BT, United Kingdom}
}

\date{April, 2023}

\begin{document}

\maketitle



\begin{abstract}
\noindent Diffuse optical tomography (DOT) is a severely ill-posed nonlinear inverse problem that seeks to estimate optical parameters from boundary measurements. In the Bayesian framework, the ill-posedness is diminished by incorporating {\em a priori} information of the optical parameters via the prior distribution.
In case the target is sparse or sharp-edged, the common choice as the prior model are non-differentiable total variation and $\ell^1$ priors. Alternatively, one can hierarchically extend the variances of a Gaussian prior to obtain differentiable sparsity promoting priors. By doing this, the variances are treated as unknowns allowing the estimation to locate the discontinuities. \\
\indent In this work, we formulate hierarchical prior models for the nonlinear DOT inverse problem using exponential, standard gamma and inverse-gamma hyperpriors. Depending on the hyperprior and the hyperparameters, the hierarchical models promote different levels of sparsity and smoothness. To compute the MAP estimates, the previously proposed alternating algorithm is adapted to work with the nonlinear model. We then propose an approach based on the cumulative distribution function of the hyperpriors to select the hyperparameters. We evaluate the performance of the hyperpriors with numerical simulations and show that the hierarchical models can improve the localization, contrast and edge sharpness of the reconstructions.
\end{abstract}

\section{Introduction}

\indent Diffuse optical tomography (DOT) is a highly ill-posed nonlinear problem which utilizes boundary measurements of near-infrared light to estimate spatially distributed absorption and reduced scattering parameters in biological tissue \cite{arridge1999optical,gibson2005recent}. Due to the highly ill-posed nature of DOT, advanced inversion techniques are required to acquire feasible estimates \cite{mozumder2014compensation,mozumder2021model,yoo2019deep}. This can be realized in a variational approach, where a regularizer based on certain assumptions is used, or through a Bayesian approach. The regularizer can, for example, utilize assumptions on smoothness or sparsity of the solution \cite{arridge1993performance,shaw2014}, or sparsity of its derivative, i.e.,  total variation (TV) \cite{paulsen1996enhanced} to obtain stable inversion. On the other hand, Bayesian estimation utilizes prior probability distributions of the unknowns, based on previously available knowledge, to compute the posterior probability distribution as a solution to the inverse problem \cite{Kaipio}. Typically the solution is computed as a single point estimate such as the \textit{maximum a posteriori} (MAP) estimate. For nonlinear problems such as the DOT, the MAP estimate lacks analytical form, therefore requiring the use of iterative methods. \\ 
\indent In Bayesian inversion, choosing a plausible prior distribution has a crucial role in mitigating the ill-posedness and determining the type of the prior information. A popular choice as prior is the Gaussian distribution, which encapsulates prior information on the mean, variance and correlations of the unknown parameters. Additionally, the Gaussian prior provides a closed-form for the posterior. In order for a Gaussian prior to provide the optimal information of a target with discontinuities, i.e., sudden jumps or sharp edges, one would need to know the locations of the discontinuities. Then, setting large variances for these locations and diminishing other variances would enable the Gaussian prior to favour discontinuities in correct locations. Evidently, knowing the locations beforehand is not the case for tomographic problems where the premise is to locate the changes in the target. \\
\indent The remedy to avoid lack of prior information and incorrect prior parameters such as the variances (of uncorrelated unknowns) is to include an uncertainty at a hierarchically higher level. By assuming the uncertain variances as unknowns, we can express the uncertainty with the hyperprior distributions. Most of the studied hierarchical models are built around the Gaussian priors, which we also investigate in this work. \\
\indent Two popular Gaussian priors used for the hierarchical models in tomographic inverse problems are the uncorrelated Gaussian prior \cite{calvetti19}, that is, a multivariate Gaussian distribution without correlations and the structural difference prior
\cite{calvetti2007gaussian}. For both priors, the hierarchical models are commonly formed by assuming the variances (for difference prior also referred to as weights) to be the unknown parameters that follow the chosen hyperprior, as in \cite{calvetti2008hypermodels,calvetti19}. Alternatively, the other parameters of Gaussian priors could also be hierarchically extended, i.e., the mean \cite{guven2005diffuse}, characteristic correlation length of Gaussian smoothness prior \cite{roininen2016hyperpriors}, or all entries of the covariance \cite{aguerrebere2017bayesian}. In this work, we limit the investigation of hierarchical models on the uncorrelated Gaussian prior and structural difference prior due to their desirable sparsity and discontinuity promoting properties as described later.\\
\indent For the two considered Gaussian priors, several different types of hyperpriors have been utilized, such as the uniform \cite{calvetti_uniform}, exponential \cite{calvetti2007gaussian}, standard gamma \cite{calvetti19} and inverse-gamma \cite{calvetti2008hypermodels} distributions. The exponential and uniform distributions can be considered non-informative hyperpriors, not incorporating (almost) any assumption on the unknowns. Whereas the standard and inverse-gamma hyperpriors can be harnessed to promote sparsity or sharp-edges \cite{calvetti2020sparse,calvetti2008hypermodels}. As it turns out, the MAP estimates with the standard and inverse-gamma hyperpriors are related to sparsity regularization, i.e., a TV or $\ell^1$ regularized problem. Conveniently, using the hyperprior approach, regularized solutions that are conventionally computed from non-smooth problems can now be computed approximately from a differentiable problem \cite{calvetti19}. The differentiability arises from the fact that the MAP estimation problems from the hierarchical models lie in the differentiable $\ell^2$ regularized framework. \\
\indent In tomographic inverse problems the hierarchical models have been used for instance, in electrical impedance tomography (EIT) \cite{liu2018image,bardsley2015randomize}, magnetoencephalography (MEG) \cite{calvetti2015_krylov_MEG,nummenmaa2007hierarchical,mattout2006meg} and linearized DOT \cite{miyamoto2011phase,shimokawa2012hierarchical,Takeaki2013,yamashita2016multi,guven2005diffuse}. A lot of the previous work concentrates on using either the structural difference prior or the uncorrelated Gaussian prior, combined with the standard gamma hyperpriors for the variances aiming to enhance the spatial accuracy of the estimates \cite{calvetti_MEG,nummenmaa2007hierarchical,shimokawa2012hierarchical}. The majority of the previously studied hierarchical models applied in tomographic problems have focused on using linear or linearized forward models, such as the Rytov approximation to linearize the forward model in DOT, see, e.g., \cite{shimokawa2012hierarchical,yamashita2016multi,Takeaki2013}. Although some work with nonlinear forward models can be found. For instance, in \cite{roininen2016hyperpriors}, the authors demonstrated the use of hierarchical Whittle-Mat\'{ern} priors or in \cite{hiltunen2009combined} the mixture Gaussian prior model with hyperpriors was applied in a nonlinear DOT problem. \\
\indent In this work, we study the effect of the exponential, standard gamma, and inverse-gamma hyperprior hierarchical models on the nonlinear two-dimensional (2D) frequency-domain DOT problem. The hierarchical models are built around the uncorrelated Gaussian prior and difference prior, for which we assume the unknown variances follow hyperpriors. The hierarchical models are evaluated with simulated piecewise constant targets that would require finding the discontinuities for optimal results, which we try to obtain via the unknown variances. The role of the hyperparameters is discussed, and a simple selection method based on the cumulative distribution function of the hyperpriors is considered. This is in contrast to previous work \cite{calvetti19,calvetti2020sparse}, where a rigorous automatized hyperparameter selection method was proposed, but only for problems with a linear model. We study the MAP estimates with different hyperparameters, which are computed from the hierarchical models by adapting the previously proposed iterative alternating sequence (IAS) algorithm \cite{calvetti2008hypermodels,calvetti2020sparse} that is modified here to work with the nonlinear forward operators. 
Additionally, we give new empirical insight into the convergence of the nonlinear IAS algorithm.\\ 
\indent The remaining sections of the paper are organized as follows. In Section \ref{sec:bayes_hierarchical}, we formulate the inverse problem in the Bayesian regime and introduce the used priors and hyperpriors. Then we formulate the nonlinear IAS algorithm for each of the hierarchical models. In Section \ref{sec:DOT}, we formulate the 2D DOT inverse problem. In Section \ref{sec:Experiments}, numerical implementations are described. In Section \ref{sec:results}, we provide the numerical results of the hierarchical models from the simulated DOT problems and discuss the empirical convergence of the nonlinear IAS. In Section \ref{sec:disc_conc}, conclusions are given.

\section{Hierarchical Bayesian models with nonlinear forward model}\label{sec:bayes_hierarchical}
Let us consider a discrete observation model of the form
\begin{equation}\label{obs_model}
y = A(x)\:+\:e,
\end{equation}
where $y\in \mathbb{R}^m$ is a discrete set of measurements, $x\in\mathbb{R}^n$ are the parameters to be estimated, $A:\mathbb{R}^n\to\mathbb{R}^m$ is a nonlinear forward operator mapping the parameters $x$ to the data space, and $e\in\mathbb{R}^m$ is the additive measurement noise. The inverse problem is then to recover the set of unknown parameters $x$ from the measured noisy data $y$.\\
\indent In the Bayesian approach to inverse problems, all the parameters are considered as random variables, and the uncertainties of their values are encoded into probability distribution models \cite{Kaipio}. The solution to the inverse problem is the posterior probability distribution which is obtained through Bayes' theorem and can be written as
\begin{equation}\label{eq:bayes_the}
    \pi(x|y) = \frac{\pi(y|x)\pi_{\rm prior}(x)}{\pi(y)} \propto \pi(y|x)\pi_{\rm prior}(x),
\end{equation}
where $\pi(y|x)$ is the likelihood density and $\pi(x)$ is the prior density. Since the measurements are known, we will neglect the normalization factor $\pi(y)$. As diffuse tomographic problems such as DOT and EIT appear as severely ill-posed, choosing appropriate prior distribution plays a principal part in overcoming the ill-posedness.\\ 
\indent For inverse problems with expensive forward operators as well as high dimensional problems, computing the entire posterior is computationally infeasible. Therefore it is common to compute a point estimate from the posterior. A commonly used point estimate is the MAP estimate, which is computed as 
\begin{equation*}
    x_{\rm MAP} = \arg \max _{x}\{\pi(x|y)\}.
\end{equation*}
To formulate the likelihood, assume mutually independent unknowns $x$ and noise $e$ that is Gaussian distributed 
\begin{equation}\label{eq:gauss}
e \sim \mathcal{N}(\mu_e,C_e), \nonumber 
\end{equation}
where $\,\mu_e\in\mathbb{R}^{m}$ is the mean, and $C_{e}\in\mathbb{R}^{m\times m}$ is the covariance matrix. The likelihood is given by \cite{Kaipio}
\begin{equation}\label{eq:likelihood}
     \pi(y|x) = \pi_{e}(y-A(x)) = \frac{1}{\sqrt{(2\pi)^m|C_{e}|}}\mathrm{exp}\left(-\frac{1}{2}\norm{L_e(y - A(x)-\mu_e)}_2^2\right),
\end{equation}
where 
$|C_{e}|$ is the determinant of $C_e$, and $L_e$  is given through the Cholesky decomposition $C_e^{-1}=L_e^{\rm T} L_e$. To simplify the notations we assume here origin centered noise ($\mu_e=0$). The MAP estimate is now the minimizer of the functional
\begin{equation}\label{eq:map_with_prior}
F(x) =  \frac{1}{2}\left\Vert L_{e}(y-A(x)) \right\Vert_2^2 - \log(\pi_{ \rm prior}(x)).
\end{equation}
This problem resembles a variational approach to regularization, where the regularization term is the negative logarithm of the prior density \cite{calvetti_bayes_overview}. \\
\indent For the choice of prior in Eq. (\ref{eq:bayes_the}) we consider two different Gaussian priors which we later extend hierarchically. These prior types have already been established to produce computationally feasible MAP estimates with the hierarchical posterior model with linear forward operators \cite{calvetti2008hypermodels,calvetti2020sparse}.\\ 
\indent Let us first consider a situation where the entries $x_j$ of the unknown vector $x$ are mutually independent and Gaussian distributed, 
\begin{equation*}
    x_j \sim \mathcal{N}\left(\mu,\sqrt{\theta_j}\right),
\end{equation*}
where $\mu\in\mathbb{R}$ is the mean and $\theta_j\in\mathbb{R}_+$ the variance. Note, that all of the unknowns $x_j$ are assumed to have same mean. Then an uncorrelated Gaussian prior is given by 
\begin{equation}\label{eq:white_noise_extended}
    \pi_{\rm prior}(x) =\frac{1}{(2 \pi)^{n/2} \sum_{j=1}^{{n}}\sqrt{\theta_{j}}} \mathrm{exp}\left(-\frac{1}{2}\sum_{j=1}^n\frac{(x_i-\mu)^2}{\theta_i}\right).
\end{equation}
Alternatively, we introduce also a difference prior of form
\begin{equation}\label{eq:smoothness}
    \pi_{\rm prior}(x) = C(\theta)\mathrm{exp}\bigg(-\frac{1}{2} \sum_{t\in D}\frac{d_t}{\theta_t}^2\bigg) =  C(\theta)\exp \left(-\frac{1}{2}\left\|\Lambda^{1 / 2} B x\right\|_2^{2}\right),
\end{equation}
where $C$ is the normalization factor depending on the variances, $D\in\mathbb{R}^{q}$ is a set that contains $q$ indice pairs $(i,j)$ determining the differences $d_{({i,j})}=x_i-x_j$, $\Lambda$ is a diagonal matrix $\Lambda=\mathrm{diag}(d_{(1,0)},\dots,d_{(n,n-1)})$ where $x_0$ is defined as zero and 
$B\in\mathbb{R}^{q \times n}$ is the difference matrix defining the differences $d$ between the adjancent unknowns, such that, $d = Bx$, $d=(d_1,\dots,d_q)$.\\
\indent For the uncorrelated Gaussian prior, the MAP estimate is the minimizer of the functional
\begin{equation}\label{eq:map_white_noise}
F(x) =  \frac{1}{2}\left\Vert L_{e}(y-A(x)) \right\Vert_2^2  + \frac{1}{2}\sum_{j=1}^n\frac{(x_j-\mu)^2}{\theta_j}.
\end{equation}
Correspondingly, for the difference prior, the MAP estimate is obtained as a minimizer of
\begin{equation}\label{eq:map_smoothness}
F(x) = \frac{1}{2}\left\Vert L_{e}(y-A(x)) \right\Vert_2^2  +\frac{1}{2}\sum_{t\in D}\frac{d_j^2}{\theta_t}.
\end{equation}
Note, that for fixed variances $\theta_j$ the normalization terms of the uncorrelated Gaussian prior \eqref{eq:white_noise_extended} and difference prior \eqref{eq:smoothness} can be omitted from the minimizations \eqref{eq:map_white_noise} and \eqref{eq:map_smoothness}. As it can be seen, both prior models yield a quadratic penalty term, the first \eqref{eq:map_white_noise} for the unknowns $x$ and the latter \eqref{eq:map_smoothness} for the differences $d$. The quadratic penalty tends to smooth sudden jumps, that is, the outliers or the edges of a piecewise constant target. For imaging modalities, this smoothness on $x$ can appear as blurriness and loss of contrast in the images, possibly making the quality of the image undesired. 

\subsection{Extending the priors hierarchically}\label{sec:hierarhcical_models}
To overcome the incomplete knowledge of the discontinuity locations, we assume the variances as random variables and aim to estimate these locations. We extend Bayes' theorem in Eq. (\ref{eq:bayes_the}), by introducing a conditional probability on the variances 
\begin{equation*}
    \pi(x|y,\theta) \propto \pi(y|x)\pi_{\rm prior}(x|\theta)\pi_{\rm hyper}(\theta),
\end{equation*}
where $\pi_{\rm hyper}(\theta)$ is now a hyperprior of the variances. Therefore the MAP estimate of the hierarchically extended posterior is given by
\begin{equation}\label{map_generic_hierarhcical}
    (x_{\rm MAP},\theta_{\rm MAP}) = \arg \max _{x,\theta}\{\pi(x|y,\theta)\}.
\end{equation}
Choosing a hyperprior for the variances allows us to express our assumption on how rare or pronounced the discontinuities are. On the other hand, assuming a uniform hyperprior provides no information on how many or extreme the discontinuities are. 

\subsubsection{Exponential hyperprior}
In this work, we only consider the case of independent variances leading to independent hyperpriors. Letting the variances be mutually independent means that the discontinuities are assumed to be sudden. For the uncorrelated Gaussian prior (\ref{eq:white_noise_extended}) the suddenness means promoting point-like outliers. Whereas the difference prior (\ref{eq:smoothness}) promotes solutions with sudden jumps, i.e., sharp-edges. \\
\indent To begin with, we describe a non-informative hyperprior for cases where we do not want to limit the amount of large variances, that is, the amount of outliers. Let us consider the exponential hyperprior for variance $\theta_j$
\begin{equation}\label{exponential}
    \pi_{\rm hyper}(\theta_j) \propto \exp \left(-\frac{1}{2}\frac{\gamma}{\theta_j}\right),
\end{equation}
where $\gamma>0$ is a hyperparameter called the rate parameter. 
If variances $\theta_j$ of the uncorrelated Gaussian prior (\ref{eq:white_noise_extended}) follow the exponential hyperprior with the same hyperparameter $\gamma$, computing the MAP estimate is equal to minimizing a functional
\begin{equation}\label{eq:map_exp}
    F(x, \theta) = \frac{1}{2}\left\|L_e(y-A(x))\right\|_2^{2}+\frac{1}{2}\sum_{j=1}^{n}\frac{(x_j-\mu)^2}{\theta_j}+\frac{\gamma}{2}\sum_{j=1}^{n}\frac{1}{\theta_j}-\frac{1}{2}\sum_{j=1}^{n}\log\theta_j.
\end{equation}
Notice, that the last term of Eq. \eqref{eq:map_exp} emerges from the normalization factor of the prior (\ref{eq:white_noise_extended}) which cannot be omitted due to variances $\theta$ being part of the estimation. Assuming the difference prior \eqref{eq:smoothness} and variances with exponential hyperprior leads minimizing 
\begin{equation}\label{eq:map_exp_difference}
    F(x, \theta) = \frac{1}{2}\left\|L_e(y-A(x))\right\|_2^{2}+\frac{1}{2}\sum_{t\in D}\frac{d_t^2}{\theta_t}+\frac{\gamma}{2}\sum_{t\in D}\frac{1}{\theta_t}-\frac{1}{2}\sum_{t\in D}\log\theta_t.
\end{equation}
\indent While increasing the hyperparameter $\gamma$, the exponential hyperprior (\ref{exponential}) starts favouring large variances over small ones, hence allowing more outliers or sharp edges to occur. On the other hand, when $\gamma \rightarrow 0$, the exponential hyperprior approaches a uniform distribution. Thus, with small $\gamma$ we can get a behaviour close to a uniform hyperprior and favour small and large variances almost equally. Note, that if $\gamma=0$, the hyperprior is uniform and the minimizer of \eqref{eq:map_exp} corresponds to the mean $x=\mu$ (or $d=0$ for \eqref{eq:map_exp_difference}). This can be observed from \eqref{eq:map_exp} and \eqref{eq:map_exp_difference}, by setting $x_j=\mu$ (or $d_j = 0$) for all $j$ and taking limit $\theta\rightarrow0$, which has no lower bound. As described later in Section \ref{sec:IAS}, the beneficial aspect of using the exponential hyperprior over the uniform hyperprior is its inherent way of incorporating the (logarithmic) positivity constraint for the variances.  

\subsubsection{Standard gamma and inverse-gamma hyperpriors}
The exponential hyperprior is easy to apply due to it having only a single hyperparameter. If the target is sparse, a possible drawback of using a non-informative hyperprior, such as the exponential (with small $\gamma)$, is the tendency to promote many outliers leading to noisy background. In order to promote sparse solutions, we need to consider hyperprior densities decreasing towards the large variances, hence favouring less outliers. As discussed in \cite{calvetti2008hypermodels}, if we wish to have only few prominent outliers, the established hyperpriors are the standard gamma distribution
\begin{equation}\label{eq:gamma_dist}
    \pi_{\rm hyper}(\theta_j|\beta,\vartheta)\propto \theta_j^{\beta-1} \mathrm{exp}\Big(-\frac{\theta_j}{\vartheta}\Big),
\end{equation}
and the inverse-gamma distribution
\begin{equation}\label{eq:inv_gamma}
   \pi_{\rm hyper}(\theta_j|\beta,\vartheta) \propto \theta_j^{-\beta-1} \exp \left(-\frac{\vartheta}{\theta_j}\right),
\end{equation}
where $\beta>0$ and $\vartheta>0$ are hyperparameters called shape and scaling parameter, respectively. As demonstrated in \cite{calvetti2008hypermodels}, the samples drawn from the standard gamma distribution are more likely to be outliers than samples drawn from the inverse-gamma distribution. On the other hand, the outliers drawn from the standard gamma distribution are less extreme compared to the inverse-gamma distribution. This indicates that the inverse-gamma hyperprior is ideal for promoting sparse solutions with pronounced outliers. In contrast, the standard gamma hyperprior works best with less sparse solutions with only moderate outliers.\\
\indent In order to understand the effect of these hyperpriors, let us consider the uncorrelated Gaussian prior (\ref{eq:white_noise_extended}), with unknown variances $\theta_j$ following standard gamma hyperprior \eqref{eq:gamma_dist}
\begin{equation*}
    \theta_j\sim\operatorname{Gamma}\left(\beta, \vartheta_{j}\right).
\end{equation*}
The MAP estimate is then the minimizer of a functional
\begin{equation}\label{eq:map_standard}
    F(x,\theta) = \frac{1}{2}\norm{L_e(y-A(x))}_2^{2}+\frac{1}{2} \sum_{j=1}^{n} \frac{(x_j-\mu)^{2}}{\theta_{j}}+\sum_{j=1}^{n}\left[\frac{\theta_{j}}{\vartheta_j}-\eta \log \left(\frac{\theta_{j}}{\vartheta_j}\right)\right],
\end{equation}
where $\eta=\beta-\frac{3}{2}$. In the previous work \cite{calvetti2020sparse}, functional \eqref{eq:map_standard} was reformulated solely in terms of $x$, i.e., by finding the function $f$ so that $\theta = f(x)$. Further, it was shown in \cite{calvetti2015_krylov_MEG} that for fixed $x$ and $\eta>0$ the functional has the limit
\begin{equation}\label{eq:gamma_lim_lemma}
    \lim_{\eta\to0}F(x,f(x)) = \frac{1}{2}\norm{L_e(y-A(x))}_2^{2} + \sqrt{2} \sum_{j=1}^{n}\frac{|x_j-\mu|}{\sqrt{\vartheta_j}},
\end{equation}
that is the $\ell^1$ constraint on $x$. In other words, by using small $\eta$, the MAP estimate approximately yields the solution of the $\ell^1$ constrained problem. The hyperparameters $\vartheta_j$ weight the penalty term, similar to the sensitivity-weighting, used to compensate sensitivity differences of the forward model (see e.g. \cite{sens_2}). Computing the non-smooth $\ell^1$ constrained solution directly from the right hand side of (\ref{eq:gamma_lim_lemma}) requires computationally expensive methods (see e.g. \cite{L1_methods}). The benefit of acquiring an approximation of the $\ell^1$ constrained solution via the functional (\ref{eq:map_standard}) is that now the term emerging from the prior and the hyperprior (the last sum of Eq. (\ref{eq:map_standard})) is differentiable. \\
\indent In case the variances of the uncorrelated Gaussian prior (\ref{eq:white_noise_extended}) follow the inverse-gamma hyperprior \eqref{eq:inv_gamma},
\begin{equation*}
\theta_j\sim \operatorname{InvGamma}\left(\beta, \vartheta_{j}\right),
\end{equation*}
then the minimized functional becomes 
\begin{equation}\label{eq:inv_gamma_map}
  F(x,\theta) = \frac{1}{2}\left\|L_e(y-A(x))\right\|^{2}_2+\frac{1}{2} \sum_{j=1}^{n}\frac{(x_j-\mu)^2}{\theta_j}+\sum_{j=1}^{n}\left[\frac{\vartheta_j}{\theta_{j}}+(\beta+\frac{3}{2}) \log \Big(\frac{\theta_j}{\vartheta_j}\Big)\right].
\end{equation}
Similarly as with Eq. \eqref{eq:map_standard}, we can write the functional \eqref{eq:inv_gamma_map} merely in terms of $x$ \cite{calvetti2020sparse}
\begin{equation}\label{eq:inv_gamma_penalty}
    F(x,f(x)) = \frac{1}{2}\left\|L_e(y-A(x))\right\|^{2}_2+\sum_{j=1}^{n} \log \left(\frac{(x_{j}-\mu)^{2}}{\vartheta_j}+2\right)^{\beta+\frac{3}{2}}.
\end{equation}
As pointed out in \cite{calvetti2020sparse}, when $\vartheta_j=1$ for all $j$ and $\beta\rightarrow0$, the logarithmic penalty term in (\ref{eq:inv_gamma_penalty}) approaches a penalty equivalent to one produced by $n$ invidual student distributions   
\begin{equation*}
    \text{Student}(x_j \mid \nu) \propto \frac{1}{\left(1+\frac{(x_j-\mu)^{2}}{\nu}\right)^{(\nu+1) / 2}},
\end{equation*}
when $\nu=2$. That is, a distribution favouring outliers. We note that, while using a linear forward operator $A$ and assuming the minimizer to be unique, the MAP minimization problem (\ref{eq:map_standard}) with the standard hyperpriors was shown to be globally convex \cite{calvetti2015_krylov_MEG}. Similarly, by observing the positive definiteness of the Hessian, the MAP minimization problem (\ref{eq:map_exp}) with the exponential hyperprior could also be shown to be strictly convex. Whereas, with the inverse-gamma hyperprior, the minimization problem (\ref{eq:inv_gamma_map}) is only locally convex, where the convexity radius depends on the hyperparameter $\beta$ \cite{calvetti2020sparse}. \\
\indent With the difference prior (\ref{eq:smoothness}), the standard gamma and inverse-gamma hyperpriors have similar effect as in Eqs. \eqref{eq:gamma_lim_lemma} and \eqref{eq:inv_gamma_penalty}, but now in terms of the differences $d$. If the variances of the difference prior follow the standard gamma hyperpriors, the MAP estimate has a limit
\begin{equation}\label{eq:gamma_lim_lemma_2}
    \lim_{\eta\to0}F(x,f(x)) = \frac{1}{2}\norm{L_e(y-A(x))}_2^{2} + \sqrt{2} \sum_{t\in D}\frac{|d_t|}{\sqrt{\vartheta_t}}.
\end{equation}
Therefore, for small $\eta$ the difference prior with variances following standard gamma hyperpriors yields a penalty term equivalent to the weighted TV penalty. Correspondingly to Eq. \eqref{eq:gamma_lim_lemma}, for the inverse-gamma hyperpriors, the minimized functional is \cite{calvetti2008hypermodels}
\begin{equation}\label{eq:inv_gamma_penalty_2}
    F(x,f(x)) = \frac{1}{2}\left\|L_e(y-A(x))\right\|^{2}_2+\sum_{t\in D} \log \left(\frac{d_t^{2}}{\vartheta_t}+2\right)^{\beta+\frac{3}{2}}.
\end{equation}
The penalty term in Eq. (\ref{eq:inv_gamma_penalty_2}) is similar to the Perona-Malik functional, used for edge weighted diffusion in image processing \cite{Perona}. By applying the standard and the inverse-gamma hyperpriors with the difference prior, we get two slightly different optimization problems, trying to establish the edges via the variances.

\subsection{Iterative alternating sequence}\label{sec:IAS}
For the efficient computation of the MAP estimates from models with many additional unknowns, such as the variances, it is common to consider alternating algorithms. The alternating algorithms tackle the minimization problem by alternatingly updating the actual unknowns and the model (prior) parameters. The framework of alternating algorithms for hierarchical models in tomography problems is well established, with various different variations (see e.g \cite{aguerrebere2017bayesian,zhang2012map,liu2018image,hiltunen2009combined,guven2005hierarchical}). In this paper, we extend the iterative alternating algorithm (IAS) described in \cite{calvetti19,calvetti2007gaussian}, to minimize the functionals (\ref{eq:map_exp}),(\ref{eq:map_standard}) and (\ref{eq:inv_gamma_map}), as well as their corresponding form for the difference prior, with a nonlinear forward operator. To our knowledge, this is the first time this algorithm has been modified for a nonlinear forward model.\\
\indent First, consider the uncorrelated Gaussian prior (\ref{eq:white_noise_extended}) and inspect the corresponding energy functional in two parts as 
\begin{equation}\label{eq:IAS_wh}
        F(x,\theta) = \left(\mathrlap{\overbrace{\phantom{\frac{1}{2}\|L_e(y-A(x))\|^{2}_2+\frac{1}{2} \sum_{j=1}^{n} \frac{\delta (x_{j}-\mu)^{2}}{\theta_{j}}}}^{\mathrm{a)}}}\frac{1}{2}\|L_e(y-A(x))\|^{2}_2+\underbrace{\frac{1}{2} \sum_{j=1}^{n} \frac{(x_j-\mu)^{2}}{\theta_{j}}+\frac{1}{2}\sum_{j=1}^{n}\log\theta_j+\sum_{j=1}^{n}g_j(\theta_j)}_{\mathrm{b)}}\right),
\end{equation}

\noindent where the functions $g_j$ depend on the selected hyperprior type. Now the a) part contains the terms depending on $x$ and the b) part depends on the variances $\theta$. For fixed $\theta$, minimizing Eq. \eqref{eq:IAS_wh} corresponds to minimizing the part a) 
\begin{equation}\label{eq:opt_1}
     \widehat{x} = \arg \underset{x}{\min}\, F(x,\theta) = \arg \underset{x}{\min}\,\frac{1}{2}\|L_e(y-A(x))\|^{2}_2+\frac{1}{2} \sum_{j=1}^{n} \frac{(x_{j}-\mu)^{2}}{\theta_{j}}.
\end{equation}
Correspondingly, for fixed $x$ we only need to minimize the part b). 
The IAS algorithm for nonlinear forward operators is outlined in Algorithm \ref{alg:nonlinear_IAS}.
\begin{algorithm}[t!]
\caption{Generic iterative alternating sequential (Generic IAS)}
\label{alg:nonlinear_IAS}
\begin{algorithmic}
    \STATE Determine hyperparameters
    \STATE Set initial $(x^0,\theta^0$) and $t=0$
    \STATE \textbf{while} solution $(x^t,\theta^t)$ not converged
    \bindent
    \STATE For fixed $\theta^t$ compute minimizer $\widehat{x}$ of  (\ref{eq:IAS_wh}) and set $x^{t+1} = \widehat{x}$
    \eindent
    \bindent
    \STATE For fixed $x^{t+1}$ compute minimizer $\widehat{\theta}$ of (\ref{eq:IAS_wh}) and set $\theta^{t+1} = \widehat{\theta}$
    \STATE $t = t + 1$
    \eindent
\end{algorithmic}
\end{algorithm}
On each iteration of the IAS, the variances are updated based on the results of the previous iteration. This idea is closely related to the bootstrap priors, where the prior is updated based on the previous recontructions \cite{bootstrap}. \\
\indent For the considered hyperpriors, the IAS is efficient since by applying it with the exponential, standard gamma or inverse-gamma hyperpriors, there exist analytical minimizer with respect to $\theta$. For instance, consider the variances following standard gamma hyperpriors. Then, for fixed $x$, minimizing functional (\ref{eq:map_standard}) equals minimizing b)
\begin{equation}\label{eq:IAS_objective}
    \widehat{\theta} = \arg \underset{\theta}{\min}\, F(x,\theta) = \arg \underset{\theta}{\min}\,\frac{1}{2} \sum_{j=1}^{n} \frac{(x_j-\mu)^{2}}{\theta_{j}}+\sum_{j=1}^{n}\left[\frac{\theta_{j}}{\vartheta_j}-\alpha \log \left(\frac{\theta_{j}}{\vartheta_j}\right)\right]
\end{equation}
which has the analytical minimizer \cite{calvetti2008hypermodels}
\begin{equation}\label{eq:stand_gamma_analytical}
    \widehat{\theta}_j = \vartheta_j\left(\frac{\eta}{2}+\sqrt{\frac{\eta^2}{4}+\frac{(x_j-\mu)^2}{2\vartheta_j}}\right).
\end{equation}
Correspondingly, for the inverse-gamma hyperpriors, with fixed $x$ the functional (\ref{eq:inv_gamma_penalty}) has the minimizer  
 \begin{equation}\label{eq:invgamma_analytical}
    \widehat{\theta}_{j}=\frac{1}{\beta+3 / 2}\left(\vartheta_{j}+\frac{1}{2} (x_{j}-\mu)^{2}\right).
\end{equation} 
Moreover, the functional (\ref{eq:map_exp}) with exponential hyperpriors  also has an analytical minimizer with respect to the variances, that is \cite{calvetti2007gaussian}
\begin{equation}\label{eq:exp_analytical}
    \widehat{\theta}_{j}=(x_{j}-\mu)^{2}+\gamma.
\end{equation} 
As it can be seen from equations (\ref{eq:stand_gamma_analytical})-(\ref{eq:exp_analytical}), for each hyperprior type, the hyperparameters determine (positive) lower bounds for the variances. For the standard gamma hyperprior \eqref{eq:stand_gamma_analytical} the shape parameter $\eta$ provides a small relaxation parameter so that the variances will not vanish, and $\vartheta_j$ controls how small the relaxation parameter is and how much the term $(x-\mu)^2$ is amplified. Whereas, for the inverse gamma hyperprior (\ref{eq:invgamma_analytical}), the scale parameter $\vartheta_j$ provides a relaxation parameter, and controlling $\beta$ can be used to amplify or dampen the variance values. The exponential hyperprior \eqref{eq:exp_analytical} has only $\gamma$ hyperparameter, which acts as the relaxation parameter while updating the variances.\\
\indent In order to apply the IAS algorithm to the MAP estimation problem with the difference prior (\ref{eq:map_smoothness}), we would first need to determine the form of the normalization factor $C$. Unfortunately, the closed-form of the normalization factor remains unknown in terms of $x$. Therefore, one option is to use Bayesian variational methods to approximate the posterior distribution, which is already established for the hierarchical models \cite{yamashita2016multi,Takeaki2013,shimokawa2012hierarchical}. Alternatively, as described in \cite{calvetti2007gaussian}, to avoid cumbersome estimation of the normalization factor, we reformulate the minimization problem in terms of the differences as
\begin{equation}\label{eq:ias_smoothness}
        F(d,\theta) = \left(\mathrlap{\overbrace{\phantom{\frac{1}{2}\|L_e(y-A(Pd))\|^{2}_2+\frac{1}{2} \sum_{j=1}^{n} \frac{\delta d^{2}}{\theta_{j}}}}^{\mathrm{a)}}}\frac{1}{2}\|L_e(y-A(Pd))\|_2^{2}+\underbrace{\frac{1}{2} \sum_{j=1}^{n} \frac{d^{2}}{\theta_{j}}+\frac{1}{2}\sum_{j=1}^{n}\log\theta_j+\sum_{j=1}^{n}g_j(\theta_j)}_{\mathrm{b)}}\right),
\end{equation}
where $P\in\mathbb{R}^{q\times q}$ is a matrix such that $x=Pd$. Then we can use Algorithm \ref{alg:nonlinear_IAS} directly to minimize (\ref{eq:ias_smoothness}). Due to fact that the difference matrix $B$ in (\ref{eq:smoothness}) is not invertible 
the matrix $P$ cannot be formulated directly from the relation $d=Bx$. Instead, the inverse mapping must be inferred from the fact that summing differences $d_i$ of any closed loop in the mesh needs to be zero. This constraint can be expressed in a matrix form as
\begin{equation}\label{eq:constraint}
    Md + \epsilon= 0,
\end{equation}
where $M\in\mathbb{R}^{p\times q}$ is a matrix representing each loop constraint in its rows, $p$ is the number of loops around the mesh elements, and $\epsilon>0$ is a small constant allowing us to compute approximated solutions of $Md=0$. The detailed formulation of the constraint \eqref{eq:constraint} is shown in \cite{calvetti2007gaussian}. The constraint (\ref{eq:constraint}) can be simultaneously solved with the part a) of (\ref{eq:ias_smoothness}) to make the differences $d$ satisfy the constraint. Compared to the hierarchical model with white noise prior, minimizing part $a)$ is now drastically more expensive due to a high number of unknown differences. For instance, a triangular mesh has $3(n-1)-b$ differences, where $b$ denotes the number of nodes on the hull.

\subsection{Selecting the hyperparameters}\label{sec:selection}
Selecting the hyperparameters is especially important for the more informative hyperpriors, which aim to favor certain solution types such as the sparsity promoting standard and inverse-gamma hyperpriors. We will  focus on the selection of the hyperparameter $\vartheta$ due to its easily interpreted role on controlling the significance of the constraints in (\ref{eq:map_standard}) and (\ref{eq:inv_gamma_penalty}). Let us first review a hyperparameter $\vartheta$ selection method for models with a linear forward operator as described in \cite{calvetti19}.

\subsubsection{Automatizing the scaling hyperparameters for linear forward models}\label{sec:auto_select}
In \cite{calvetti19} it was shown that for a linear forward operator $A(x)=Ax$ satisfying the exhangeability condition, the hyperparameters $\vartheta_j$ should be set as
\begin{equation}\label{eq:hyperselection}
    \vartheta_j=\frac{D}{\left\|A e_{j}\right\|^{2}_2}, 
\end{equation}
where term $D$ depends on the estimated signal-to-noise ratio, hyperparameter $\beta$, and assumed sparsity of the $x$ (see \cite{calvetti19} for details).\\
\indent As discussed in \cite{calvetti19}, the idea of setting the hyperparameters $\vartheta$ according to the forward model is closely related to the idea of sensitivity-weighting used in the regularization scheme. In sensitivity-weighting the regularizing norm such as the $\ell^1$ regularizer is scaled for each component $x_j$ with respect to the sensitivity 
\begin{equation*}
\left\|\frac{\partial(A x)}{\partial x_{j}}\right\|_2=\left\|A e_{j}\right\|_2,
\end{equation*} 
where $e_j$ is the $j$th canonical basis vector. Similarly, in the automatized selection (\ref{eq:hyperselection}), the sensitivity term $\left\|A e_{j}\right\|_2$ appears. The sensitivity-weighting is established to improve the spatial quality of recontructions in applications, such as MEG \cite{calvetti_MEG}, where the sensitivity differences of the linear forward model are prominent. On the other hand, if we choose the same value for all $\theta_j$, no spatial preferences are incorporated among the unknowns $x$.

\subsubsection{Choosing the hyperparameters based on the confidence interval}\label{sec:confidence}
For a nonlinear model, we cannot deduce similar selection method utilizing the sensitivity of the model (i.e. $\left\|A e_{j}\right\|$), since the sensitivity of $A(x)$ varies for different $x$. In practical applications it is common to know coarse upper bound for the unknown physical quantities. Therefore, we can incorporate this information via the confidence intervals of the probability distributions. For instance, let us assume that $\abs{x_j-\mu}<M$, $M>0$ with 95\% probability. Then the variance $\theta_j$ is set as $U=(M/2)^2$. \\
\indent On the other hand, if $\theta_j$ follows the hyperprior, we have 5\% probability for $x_j$ being an outlier, for which $\theta_j$ should be larger than $(M/2)^2$. Thereby we can determine $\vartheta_j$ from the cumulative distribution function (CDF) of the hyperprior, i.e., we choose $\vartheta_j$ s.t.
\begin{equation}\label{eq:cdf}
    \mathrm{CDF}(U,\vartheta_j) = 0.95. 
\end{equation}
For the exponential hyperprior, the CDF is not bounded and thus cannot be used to determine $\gamma$. 
\indent \revision{
\begin{remark}
As discussed in Section \ref{sec:auto_select}, with a linear model, the automatized selection works as a sensitivity-weighting, compensating the sensitivity differences between the estimated unknowns. We would like to remark that similar sensitivity compensating could be achieved in the nonlinear case with more ad hoc hyperparameter selection. For instance, if we know that the unknowns near the sources/detectors are more sensitive, the hyperparameters could be scaled based on spatial distances. Alternatively, one could evaluate $A(x)$ for multiple points $x_a$ and try to approximate the sensitivities
\begin{equation*}
\left\|\frac{\partial A (x_a)}{\partial x_{j}}\right\|_2=\left\|J(x_a) e_{j}\right\|_2,
\end{equation*} 
by utilizing the Jacobian $J$ at point $x_a$.
 \end{remark}
}

\section{Diffuse optical tomography}\label{sec:DOT}
Diffuse optical tomography (DOT) is a technique for imaging spatially varying optical parameters in biological tissue \cite{arridge1999optical,gibson2005recent,durduran2015}. 
The distribution of these optical parameters provides tissue biochemical and structural information with applications, for example, in early-diagnosis and imaging of breast cancer \cite{grosenick2016}, monitoring neonatal brain health \cite{hebden2002three}, functional brain imaging of adults \cite{wheelock2019high}, and preclinical imaging of small animals \cite{chen2012optical}. In general, DOT is non-ionizing and non-invasive, and its instrumentation is relatively simple, low-cost, and portable compared to conventional medical tomographic techniques \cite{wheelock2019high}. For more information on image reconstruction problem of DOT and various methodologies, see e.g. \cite{arridge2009,durduran2015,hoshi2016,shaw2014,diamond2006} and the references therein.

\subsection{Forward model}

\indent In a typical DOT measurement setup, near-infrared light is introduced into an object from its boundary. Let $\Omega \subset \mathbb{R}^d, \, (d=2 \, \rm{or} \, 3)$ denote the domain with boundary $\partial \Omega$ where $d$ is the (spatial) dimension of the domain. In a diffusive medium, like soft biological tissue, the commonly used light transport model for DOT is the diffusion approximation to the radiative transfer equation \cite{Ishimaru}. Here, we consider the frequency-domain version of the diffusion approximation \cite{arridge1999optical}
\begin{equation}\label{deeqn}
 \left(-\nabla \cdot \frac{1}{{d}(\mua(r)+\mus(r))} \nabla  + \mua(r) + \frac{{\rm j}\omega}{c} \right) \Phi(r) = 0, \quad r \in \Omega,
\end{equation}
\begin{equation}\label{deeqnbn}
\Phi(r)+\frac{1}{2\zeta}\frac{1}{d(\mua (r)+\mus(r))} \alpha \frac{\partial \Phi(r)}{\partial \hat{n}} = \left\{ \begin{array}{ll}
         \frac{q}{\zeta}, & r \in s\\
         0, & r \in \partial \Omega \setminus s \end{array} \right. ,
\end{equation}
where $\Phi(r)$ is the photon density, $\mua(r)$ is the absorption coefficient, $\mus(r)$ is the reduced scattering coefficient, j is the imaginary unit and $c$ is the speed of light in the medium. The parameter $q$ is the strength of the light source at location $s \subset \partial\Omega$, operating at angular modulation frequency $\omega$. Further, the parameter $\zeta$ is a dimension-dependent constant ($\zeta$ = $1/\pi$ when $\Omega \subset \mathbb{R}^{2}$, $\zeta$ = $1/2$ when $\Omega \subset \mathbb{R}^{3}$) and $\alpha$ is a parameter governing the internal reflection at the boundary $\partial \Omega$, and $\hat{n}$ is an outward unit vector normal to the boundary. 
The measurable data on the boundary of the object,  exitance $\Gamma(r)$, is given by
\begin{equation}\label{bdd}
\Gamma(r) = -\frac{1}{d(\mua (r)+\mus(r))} \frac{\partial \Phi(r)}{\partial \hat{n}} = \frac{2\zeta}{\alpha}\Phi(r). 
\end{equation}
\indent The numerical approximation of the forward model (\ref{deeqn})-(\ref{bdd}) is often based on a finite element (FE) approximation \cite{arridge1999optical}. In this work, we use a finite dimensional approximations in piecewise linear basis for absorption, reduced scattering and fluence as described in \cite{arridge1999optical}. 

\subsection{Inverse problem of DOT}
Let us consider observation model (1), where y is the vector of measurable data that typically in frequency domain DOT are logarithm of amplitude and phase delay of exitance, $x=(\mu_{\rm{a}_{,1}},\dots,\mu_{\rm {a}_{,n/2}},\mu'_{\rm{s}_{,1}},\dots,\mu'_{\rm{s}_{,{n/2}}})\in \mathbb{R}^{n}$ is vector containing the absorption $\att$ and reduced scattering $\scat$, $A(x)$ is the discretised forward operator, i.e. the FE-approximation of the forward model (\ref{deeqn})-(\ref{deeqnbn}), and $e$ is the additive noise. For this nonlinear observation model, we can compute the MAP estimates of the hierarchical models described in Section \ref{sec:hierarhcical_models} with the introduced nonlinear IAS algorithm (Algorithm \ref{alg:nonlinear_IAS}).\\
\indent In order to use the IAS iterations, we first need to minimize part a) of the functional (\ref{eq:IAS_wh}). Since the minimizer has no closed-form solution, iterative methods, such as Gauss-Newton method are used to approximate the solution. The Gauss-Newton iterations can be written as
\begin{equation}\label{eq:gn_update}
x^{i+1} = x^{i} + s_i\delta x^{i}
\end{equation}
with step-length parameter $s_i$. The update direction $\delta x^{i}$ is given by
\begin{equation}\label{eq:deri}
\begin{split}
     \delta x^{i}  = \left( J_A^{\rm T}C_e^{-1}J_A^{\rm T} +  C_{x}^{-1}\right)
\left(J_A^{\rm T}C_e^{-1}(y-A(x))
 + C_x^{-1} (x-\mu) \right),
\end{split}
\end{equation}
where the prior covariance 
\begin{equation*}
C_x = \begin{pmatrix} C_{\mu_{\rm a}} & 0 \\ 0 &  C_{\mu'_{\rm s}}
\end{pmatrix} 
\end{equation*}
contains separate prior covariances for the absorption $C_{\mu_{\rm a}}$ and reduced scattering $C_{\mu'_{\rm s}}$. Here Jacobian $J_A$ is the discrete representation of the Fr{\'e}chet derivative of the nonlinear operator $A(x)$ at the point $x^{i}$. For the difference prior model we need to estimate the differences $d$ by substituting $d=Px$. After this the chain rule is used to obtain the Jacobian of $A(Pd)$ with respect to the difference $d$. The following minimization step (part b) to update the variances is then performed by computing the analytical minimizers, using either (\ref{eq:stand_gamma_analytical}), (\ref{eq:invgamma_analytical}) or (\ref{eq:exp_analytical}) depending on the used type of hyperprior.

\section{Simulations}\label{sec:Experiments}

The DOT simulation domain was set to a circle with a radius of $25\,$mm. The setup consisted of $32$ sources and $32$ detectors modeled as $2\:\rm{mm}$ wide surface patches located at equispaced angular intervals on the boundary. Hence the total number of data was 2048, that is, all source-detector pairs of the logarithm of amplitude and phase. For data simulation, the domain was discretized using 10062 triangular elements and $n_s=5149$ nodes. The solution of the DA (\ref{deeqn})-(\ref{deeqnbn}) was numerically approximated using the finite element method via the the Toast++ software \cite{TOAST}. The simulated data was corrupted with additive white noise that was drawn from Gaussian distribution $\pi(e)\sim\mathcal{N}(0,C_e)$, where standard deviations ($C_e=\mathrm{diag} (\sigma_1^2,\dots,\sigma_n^2)$) of the noise were selected as 0.4\% of the maximum value from (complex and real parts) of the noise-free measurement data. For all simulated targets, the signal-to-noise ratio was approximately (in decibels) $45$.\\ 
\indent In the inverse problem, the FE-approximation of the DA implemented with Toast++ was used to approximate the model for light propagation. The FE-mesh of the photon density contained 9102 triangular elements and $n=4663$ nodes. For the image space, we used the same FE-mesh.  
Since a reduced number of nodes was used for the inversion, a small modeling error was included in the simulations.\\
\indent The inner iterations were performed by the Gauss-Newton method. A bisection style line-search (see \cite{line_search}), which we have previously found to be well suited for step length selection in DOT, was used for determining the step length of the algorithm. 
The stopping criteria for both linear and nonlinear IAS was met when the relative difference of the estimates $x$ was less than $\varepsilon$, i.e.,
\begin{equation}\label{criteria}
    \frac{\norm{x^t-x^{t+1}}_2}{\norm{x^t}_2}<\varepsilon, 
\end{equation}
over three consecutive iterations or when the minimized functional $F$ stopped decreasing, i.e.,
\begin{equation}\label{criteria2}
    F(x^t,\theta^t)-F(x^{t+1},\theta^{t+1})<0
\end{equation}
over three iterations. For all simulations, we used $\varepsilon=10^{-5}$. For the Gauss-Newton algorithm, the iterations were set to stop after the difference $F(x^t,\theta^t)-F(x^{t+1},\theta^{t+1})$ was smaller than $10^{-12}$ guaranteeing strict convergence of the inner iterations. The initial variances $\theta_j$ were set as $0.25^2$ and $0.0025^2$ for reduced scattering and absoprtion. \\
\indent In addition to visual inspection, the computed MAP estimates $x_{\rm{MAP}}\in\mathbb{R}^n$ were compared by computing relative errors as
\begin{equation*}
\operatorname{RE} = \frac{\norm{x_{\text{true}}-x_{\rm{MAP}}}_2}{\norm{x_{\text{true}}}_2},
\end{equation*}
where $x_{\text{true}}\in\mathbb{R}^{n_s}$ is the simulated true value. Since different parameter spaces of $x_{\text{true}}$ and  $x_{\rm{MAP}}\in\mathbb{R}^n$, the relative error is computed after interpolating both values to a basis with equal amount of nodes.\\
\indent The computations were performed with a laptop computer, equipped with Intel Core i7-11800H @ 2.30GHz processor and Nvidia T1200 Laptop GPU. The FE-solution was computed by utilizing the Toast++ software used via MATLAB (R2021b, Mathworks, Natick, MA).

\section{Results and discussion}\label{sec:results}

\begin{table}[b!]
\caption{Used hyperparameter ($\vartheta$ and $\gamma$) values for the hierarchical uncorrelated Gaussian prior (unc.) and structural difference prior (dif.). The hyperparameter $\vartheta$ values were computed from the CDF (\ref{eq:cdf}). The hyperparameter values for reduced scattering (scat.) were computed with $M\in\{0.3,\,1,\,10\}$ (unc.) and $M\in\{1,\,5,\,10\}$ (dif.) for standard gamma and $M\in\{0.3,\,1,\,5\}$ (unc.) and $M\in\{0.25,\,1,\,4\}$ (dif.) for inverse-gamma. For absorption (abs.), the assumed $M$ was 0.01 times the corresponding reduced scattering values.}
\label{tab:hp_values}
\makebox[\linewidth]{
\begin{tabular}{lllllll} 
\hline
\multirow{2}{*}{\diagbox{Hyperprior}{Hyperparam.}} & \multicolumn{2}{l}{\hspace{1cm} Low} & \multicolumn{2}{l}{\hspace{0.5cm} Intermediate} & \multicolumn{2}{l}{\hspace{1.1cm} High} \\ 
\cline{2-7}
 & Unc. &  Dif.  & Unc.  & Dif. & Unc.  &  Dif.                 \\ 
\hline
\multirow{2}{*}{Exponential ($\gamma)$} \qquad  Scat. &  $10^{-{10}}$  & $9\cdot10^{-4}$ &  2.5$\cdot 10^{-3}$  & 3.6$\cdot 10^{-3}$ & 0.25 &  0.14                           \\ \qquad \qquad \qquad \qquad \,\, Abs.
&  $10^{-{14}}$  & $9\cdot10^{-8}$ &  2.5$\cdot 10^{-7}$  & 3.6$\cdot 10^{-7}$ & 2.5$\cdot10^{-3}$ &  1.4$\cdot10^{-3}$  \\
\hline
 \multirow{2}{*}{Standard gamma ($\vartheta$)}     Scat. & 5.8$\cdot 10^{-3}$ & 6.4$\cdot 10^{-2}$
 &  6.4$\cdot 10^{-2}$  &  1.6 &  6.4  &  6.4                                 \\
\qquad \qquad \qquad \qquad \,\, Abs. & 5.8$\cdot 10^{-7}$ & 6.4$\cdot 10^{-6}$
 &  6.4$\cdot 10^{-6}$  &  1.6 $\cdot 10^{-4}$ &  6.4 $\cdot 10^{-4}$  &  6.4 $\cdot 10^{-4}$                                \\
 \hline
 \multirow{2}{*}{Inverse-gamma ($\vartheta$)} \,\, Scat. &   $4\cdot10^{-3}$  &  5.5$\cdot10^{-3}$  &  1.1$\cdot 10^{-2}$  &  8.8$\cdot 10^{-2}$  & 0.4  & 1.4                   \\
 \qquad \qquad \qquad \qquad \,\, Abs. &   $4\cdot10^{-7}$  &  5.5$\cdot10^{-7}$  &  1.1$\cdot 10^{-6}$  &  8.8$\cdot 10^{-6}$  & 4 $\cdot 10^{-3}$ & 1.4 $\cdot 10^{-4}$                  \\
\hline
\end{tabular}
}
\end{table}

\subsection{Reconstruction utilizing the uncorrelated Gaussian prior model}
First, we tested the performance of the uncorrelated Gaussian prior with the standard gamma, inverse-gamma, and exponential hyperpriors. As the test phantom, we used a piecewise constant target with two inclusions of different sizes, which should be an ideal target for the sparsity promoting standard gamma and inverse-gamma hyperpriors. \\
\indent First of all, we note that changing $\eta$ and $\beta$ parameters of the standard and inverse-gamma hyperpriors was not observed to have a significant effect on the reconstructions as long as the values were sufficiently small, i.e., $\eta<10^{-2}$ and $\beta<5$. For larger values, the effect of the hyperpriors diminished. Thereby we only provide results with hyperparameter values $\eta = 10^{-4}$ and $\beta = 1.5$ and alter the hyperparameters $\vartheta$ and $\gamma$ instead.\\
\indent For each of the hyperprior type, three different magnitudes of hyperparameters ($\vartheta$ and $\gamma$) were used: \textit{(i)} A \emph{low} value that provided only slight effect to the optimization problem, \textit{(ii)} an \emph{intermediate} value close to optimal, and \textit{(iii)} \emph{high} value strongly affecting on the optimization. The used hyperparameters $\vartheta$ and $\gamma$ are given in Table \ref{tab:hp_values}. The utilized $\vartheta$ values were chosen by using the cumulative distribution function (\ref{eq:cdf}) by assuming $|x_j-\mu|<M$ with $M\in\{0.3,\, 1,\, 10\}$ for standard gamma and $M\in\{0.3,\, 1,\, 5\}$ for inverse-gamma. The smallest $M$ value corresponds to the \emph{low}, the middle to the \emph{intermediate} and the largest to the \emph{high}. The mean $\mu$ of the uncorrelated Gaussian prior was set to be $0.01$ for the absorption and 1 for the reduced scattering. These values were the true values of the background, 
\revision{which were also used as the initial guess for the optical parameters. A few additional experiments were run with initial values slightly deviating from the background values. However, the effect of the deviation was observed to be negligible. Additionally, for inverse-gamma hyperprior we experimented with the strategy suggested in \cite{calvetti2020sparse} to run a few iterations with standard gamma hyperprior and then switch to inverse-gamma hyperprior. However, we observed no significant effect on the performance of the inverse-gamma hyperprior.}

\indent Since the magnitude of the absorption values was assumed to be 100 times smaller than the reduced scattering values, $M$ values used for the absorption were the reduced scattering values $M$ divided by 100. For the same reason, the hyperparameters $\gamma$ (which sets the minimum variance in IAS iterations \eqref{eq:exp_analytical}) used for the absorption were 0.01\% of $\gamma$ used for the reduced scattering.

\begin{figure}[tbp!]
\centering

\subfloat{\includegraphics[trim={0 0cm 0cm 0.62cm},clip,scale=0.715]{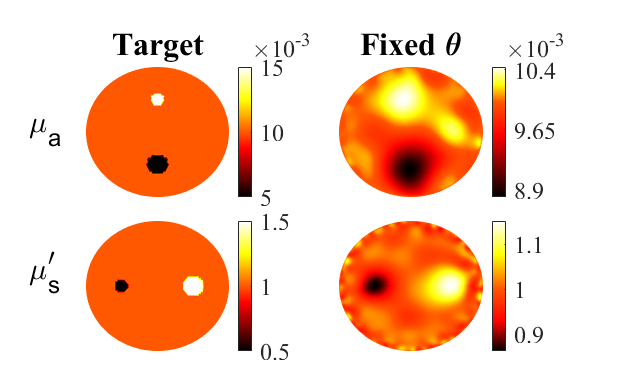}}\hspace{1cm}%

\subfloat{\includegraphics[trim={0 0cm 0.97cm 0.3cm},clip,scale=0.715]{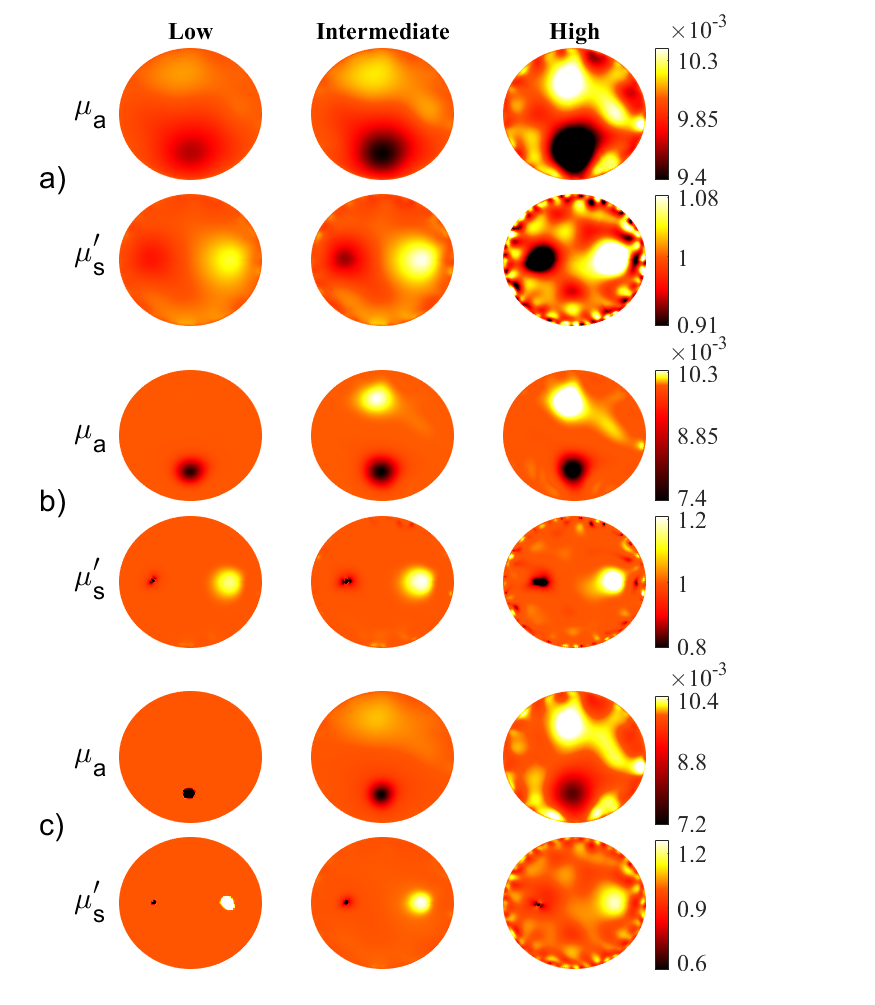}}

\caption{Computed MAP estimates. The two top-most rows show the true values and the estimates when the uncorrelated Gaussian prior (\ref{eq:white_noise_extended}) with fixed variances of $0.25^2$ and $0.0025^2$ for reduced scattering and absorption were used. Other rows show the estimates when variances were a) inverse-gamma,  b) standard gamma, and c)  exponentially distributed. The columns show the estimates with small (left), intermediate (mid) and large (right) hyperparameter $\vartheta$ or $\gamma$ values as listed in Table \ref{tab:hp_values}. The colorbars exclude some of the highest and smallest values.}
\label{fig:fig_wn}
\end{figure}

\indent Figure \ref{fig:fig_wn} shows the simulated target and the MAP estimates computed for all three cases by using the IAS method (Algorithm \ref{alg:nonlinear_IAS}). For comparison, the reconstruction was also computed using fixed variance values shown at the top of Figure \ref{fig:fig_wn}. The fixed variances were set as $0.25^2$ and $0.0025^2$ for reduced scattering and absorption, respectively, since these values produced admissible reconstruction. Note, that the reconstruction with the

\noindent fixed variances corresponds to the first iterations of the IAS. The colorbars in Figure \ref{fig:fig_wn} exclude some of the largest and smallest values. Relative errors of the estimates are shown in Table \ref{tab:wn} 

\begin{table}[b!]
\centering
\caption{Relative errors (\%) of the MAP estimates when the uncorrelated Gaussian prior (\ref{eq:white_noise_extended}) was used with low, intermediate and high hyperparameter $\vartheta$ (or $\gamma$) values. The left-most column indicates the type of the hyperprior used for the unknown variances. The non-hierarchical model with the fixed variances is denoted as "no hyperprior". For reduced scattering and absorption, the smallest relative error is bold.}
\label{tab:wn}
\makebox[\linewidth]{
\begin{tabular}{lllllll} 
\hline
\multirow{2}{*}{\diagbox{Hyperprior}{Hyperparam.}} & \multicolumn{2}{l}{Low} & \multicolumn{2}{l}{Intermediate} & \multicolumn{2}{l}{High} \\ 
\cline{2-7}
  & RE ($\mu_{\rm{a}}$)  &  RE ($\mu_{\rm{s}}^\prime$)  & RE ($\mu_{\rm{a}}$)  & RE ($\mu_{\rm{s}}^\prime$) & RE ($\mu_{\rm{a}}$)  & RE ($\mu_{\rm{s}}^\prime$)   \\ 
\hline
No hyperprior   & 7.35  &  6.80    &      7.35  &  6.80     & 7.35  &  6.80                                    \\ 

Exponential & 7.87  &  16.46 &  6.55  &  5.42    & 6.81  &  7.01                            \\ 
 Standard gamma  &   6.63  &  6.51    &  6.29  &  5.81 & \textbf{5.57}  &  \textbf{5.28}                                    \\
 Inverse-gamma  &   7.94   & 7.69   &   7.73  &  7.34  & 7.06  &  6.74
                      \\
\hline
\end{tabular}
}
\end{table}

\begin{figure}[tbp!]
\noindent\makebox[\textwidth]{\includegraphics[trim={0 0cm 1cm 0},clip,scale=0.66]{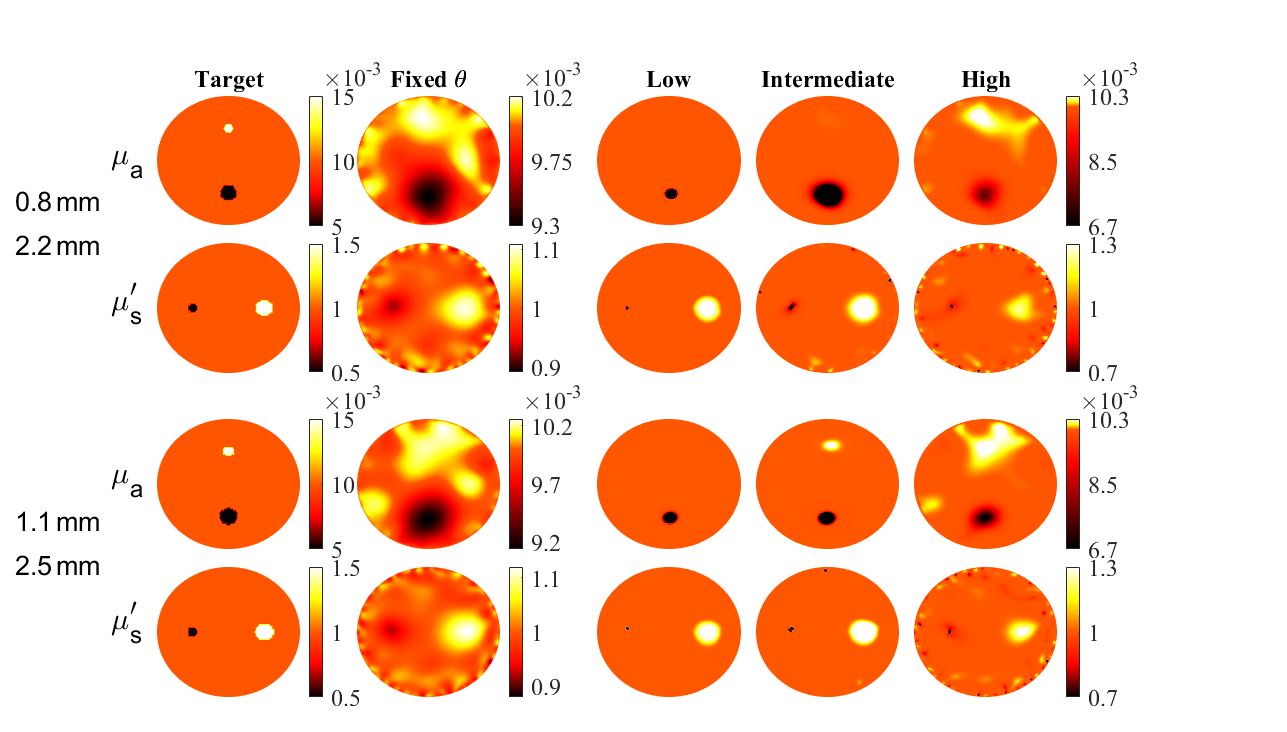}}
\caption{\revision{Reconstructions of smaller inclusions using uncorrelated Gaussian prior (\ref{eq:white_noise_extended}) with standard gamma hyperpriors. The leftmost column shows the used targets, where the radiuses of the inclusions were 2.5 mm, 1.1 mm (lower) and 2.2 mm, 0.8mm (upper). The remaining columns from left to right show the reconstructions with fixed variance values and variances with low, intermediate and high hyperparameter values. The fixed values and the hyperparameter values are the same as used in Figure \ref{fig:fig_wn} row b).}}
\label{fig:dif_radiusses}
\end{figure}

\indent From Figure \ref{fig:fig_wn}, we can observe that the reconstructions with smaller hyperparameters, i.e., with stronger prior information, produce a smooth background for both optical parameters. Whereas the reconstructions with the larger hyperparameters and weaker prior information produced inclusions with increased contrast with a noisier background. On the other hand, the strong sparsity assumption (left-most column) causes the smaller inclusions to vanish. Notably, none of the MAP estimates produced cross-talk artefacts between scattering and absorbtion, which commonly occur in DOT. For some of the reconstructions in Figure \ref{fig:fig_wn}, the single smallest reduced scattering value was drastically smaller than the others. \\
\indent Comparing the qualitative performance between the hyperprior types shows that the standard gamma hyperprior (b) yields superior absorption estimates compared to the exponential (c) and inverse-gamma (a) hyperprior, which cannot properly localize the smaller absorption inclusion. Additionally, a well-working set of hyperparameters was observed to be broad for the model with the standard gamma hyperpriors. Table \ref{tab:wn} shows that the standard gamma hyperprior also yields the best estimates regarding the relative errors.\\
\indent \revision{The uncorrelated Gaussian prior with standard or inverse-gamma hyperpriors ideally work for sparse spiky targets.  However, we would like to remark that the diffusive nature of DOT limits the size of targets that we can expect to recover, compare to \cite{isaacson1986distinguishability}. Thus, recovery of targets much smaller than those used in Figure \ref{fig:fig_wn} can not be expected. 
To demonstrate this, we tested two smaller inclusions with radiuses of 2.5 mm, 1.1 mm (larger/smaller) and 2.2 mm, 0.8 mm, whereas radiuses in Figure \ref{fig:fig_wn} were 3 mm, 1.5 mm. These targets were tested with the standard gamma hyperprior using the previously selected \textit{low}, \textit{intermediate} and \textit{high} values (Table \ref{tab:hp_values}). The reconstructions are shown in Figure \ref{fig:dif_radiusses}. From the bottom row of Figure \ref{fig:dif_radiusses} we can see that the slightly smaller inclusions (cl. Figure \ref{fig:fig_wn}) already get challenging to reconstruct. Whereas with the smallest inclusions (upper row) the reconstructions become unfeasible, either not showing the smaller inclusions or overestimating the size of the larger inclusions.}

\subsection{Reconstruction utilizing the difference prior model}
The second prior model considered in this work was the difference prior (\ref{eq:smoothness}). To test the performance of the hierarchical difference prior models, we used a target with two inclusions with sharp edges that should be ideal for investigating the edge-preserving properties. Compared to the inclusions used for the uncorrelated Gaussian prior Figure \ref{fig:fig_wn}, the inclusions were now slightly larger. This was mainly because the difference prior model was observed to vanish too small inclusions. \\
\indent Similarly, as with the uncorrelated Gaussian prior, the hierarchical difference prior was tested with three different values of hyperparameters $\vartheta$ and $\gamma$. To compute the hyperparameters $\vartheta$ from the CDF (\ref{eq:cdf}) we used three different $M$ values. For the reduced scattering we set $M\in\{1, \, 5, \, 10\}$ with standard gamma and $M\in\{0.25,\, 1, \,4\}$ with inverse-gamma hyperpriors. Since absorption values were known to be approximately $1\%$ of the reduced scattering values, the differences $d$ of absorption are approximately $1\%$ of the reduced scattering differences as well. Therefore we used 100 times smaller $M$ values for the absorption. For the non-hierarchical model, the variances were set as $0.1^2$ for reduced scattering and $0.001^2$ for absorption yielding an acceptable reconstruction quality. For $\eta$, we used $10^{-4}$ and for $\beta=3/2$, i.e., the same values as with the uncorrelated Gaussian prior. \revision{For constant $\epsilon$ needed in the compatibility condition \ref{eq:constraint}, we used value $10^{-7}$ for absorption and $10^{-5}$ for scattering. With these values, the compatibility condition was observed to be sufficiently satisfied, showing no artefacts in the reconstructions. Whereas using 20 times larger $\epsilon$ already started to cause visible artefacts in the reconstructions.} 

\begin{table}[tbp!]
\centering
\caption{Relative errors (\%) of the MAP estimates when the difference prior (\ref{eq:smoothness}) was used with low, intermediate and high hyperparameter $\vartheta$ (or $\gamma$) values. The left-most column indicates the type of the hyperprior used for the unknown variances. The non-hierarchical model with the fixed variances is denoted as "no hyperprior". For reduced scattering and absorption, the smallest relative error is bold.}
\label{tab:smoothness}
\makebox[\linewidth]{
\begin{tabular}{lllllll} 
\hline
\multirow{2}{*}{\diagbox{Hyperprior}{Hyperparam.}} & \multicolumn{2}{l}{Low} & \multicolumn{2}{l}{Intermediate} & \multicolumn{2}{l}{High} \\ 
\cline{2-7}
 & RE ($\mu_{\rm{a}}$)  &  RE ($\mu_{\rm{s}}^\prime$)  & RE ($\mu_{\rm{a}}$)  & RE ($\mu_{\rm{s}}^\prime$) & RE ($\mu_{\rm{a}}$)  & RE ($\mu_{\rm{s}}^\prime$)                 \\ 
\hline
No hyperprior   & 9.78  &  7.17   & 9.78   & 7.17   & 9.78 &   7.17                                    \\ 
Exponential &   11.14  &  9.00 &  9.60  &  6.89 & 9.44 &   6.67                           \\ 
 Standard gamma  & 10.55 & 6.54
 &  8.74  &  \textbf{4.28}  &  \textbf{7.80}  &  4.40                                 \\
 Inverse-gamma  &   11.08  &  6.14   &  10.26  &  7.67  & 9.12  &  6.34                      \\
\hline
\end{tabular}
}
\end{table}

\begin{figure}[tbp!]


  \centering
 \subfloat{\label{fig4:a}}{\includegraphics[trim={0 0cm 0cm 0.62cm},scale=0.715]{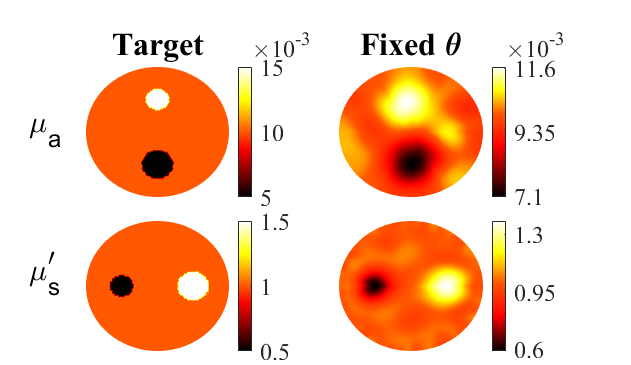}}\hspace{1em}%
  \subfloat{\label{fig4:b}}{\includegraphics[trim={0 0cm 0.97cm 0},clip,scale=0.715]{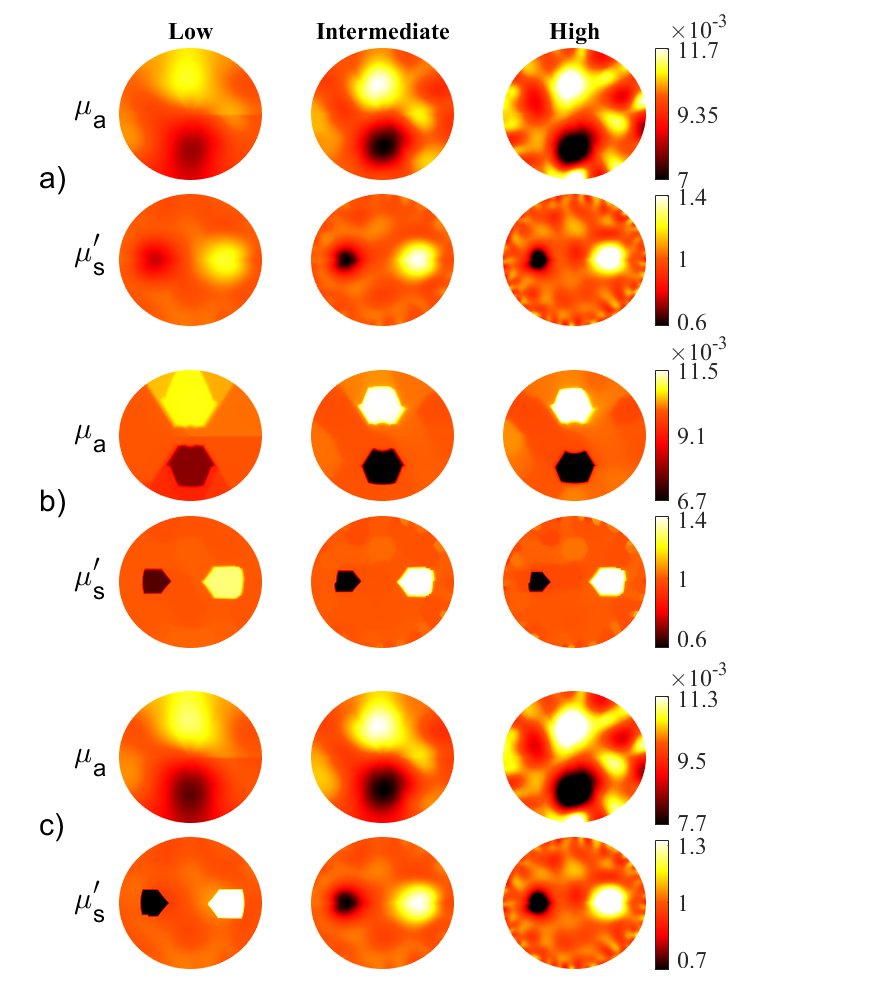}}
\caption{Computed MAP estimates. The two top-most rows show the true values and the estimates when the \revision{difference prior (\ref{eq:white_noise_extended})} with fixed variances of \revision{$0.1^2$ and $0.001^2$} for reduced scattering and absorption were used. Other rows show the estimates, when variances were a) inverse-gamma,  b) standard gamma and c)  exponentially distributed. The columns show the estimates with small (left), intermediate (mid) and large (right) hyperparameter $\vartheta$ or $\gamma$ values as listed in Table \ref{tab:hp_values}. The colorbars exclude some of the highest and smallest values.}
\label{fig:fig_smooth}

\end{figure}

\begin{figure}[t]
\noindent\makebox[\textwidth]{\includegraphics[trim={0 0cm 0 0},clip,scale=0.66]{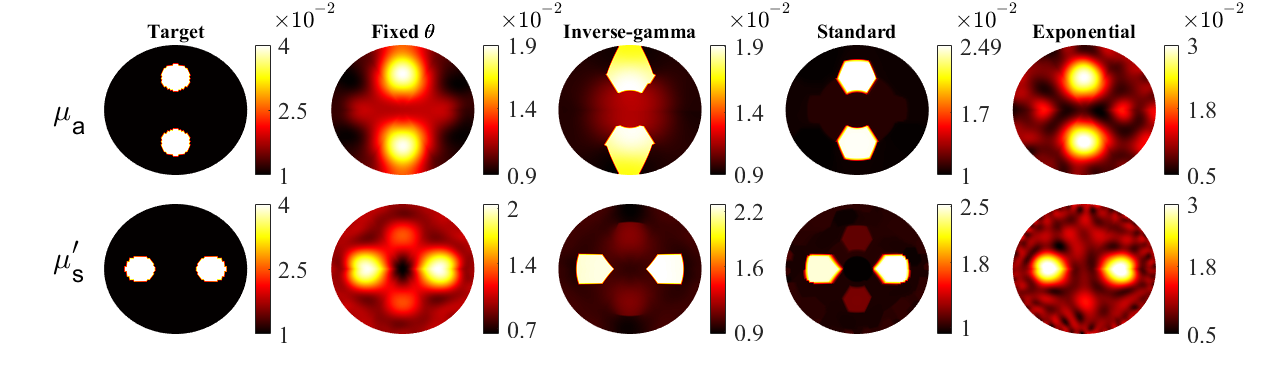}}
\caption{Computed MAP estimates of alternative target type with only positive inclusions. The left-most column show the target. The second column shows the estimates when the difference prior (\ref{eq:smoothness}) with fixed variances of $0.1^2$ and $0.001^2$ for reduced scattering and absorption was used. Other columns show the estimates, when variances were inverse-gamma, standard gamma and exponentially distributed. The hyperparameters $\vartheta$ were set as the intermediate values, also used for recontructions in Figure \ref{fig:fig_smooth}.}
\label{fig:alternative_target}
\end{figure}

\indent Figure \ref{fig:fig_smooth} shows the MAP estimates and the used target. Some of the highest and lowest values are excluded from the colorbars. Table \ref{tab:smoothness} shows the corresponding relative errors. As can be seen from the figure, the hierarchical models with the smaller hyperparameters produce sharper-edged inclusions. The disadvantage of using too small hyperparameter values is the falsely estimated edges of the absorption inclusion. With the large hyperparameters, the hierarchical models tend to increase the background noise, excluding the standard gamma hyperprior model (b), which still enhances the sharpness of the inclusions compared to the fixed variances estimates. Again the standard gamma hyperprior is observed to produce the most feasible estimates when $\vartheta$ is sufficiently large. The exponential hyperprior (c) is seemingly performing well with the reduced scattering, but the absorption inclusions remain as blurry as with the fixed variances. With the inverse-gamma hyperprior (a), no significant improvement is observed. Similarly as with the uncorrelated Gaussian noise, no cross-talk artefacts appeared in any of the reconstructions shown in Figure \ref{fig:fig_smooth}. \\
\indent Figures \ref{fig:fig_wn} and \ref{fig:fig_smooth} demonstrate the unstable performance of the exponential and inverse-gamma hyperpriors. With smaller hyperparameters, the exponential and inverse-gamma hyperpriors diminished the background noise but also significantly decreased the contrast of the inclusions. While with the larger hyperparameters, the contrast of the inclusions was better, but now the background noise was more substantial. On the other hand, the effect of these hyperpriors can be more plausible with a different kind of test target. For instance, in Figure \ref{fig:alternative_target}, we show one alternative target, with only symmetric and positive inclusions with eight times higher contrast. The MAP estimates shown in Figure \ref{fig:alternative_target} were computed using the hierarchical difference prior with the same intermediate hyperprior values, used to obtain the MAP estimates in the middle column of Figure \ref{fig:fig_smooth}. The estimates with the fixed variances were computed with variances of $0.5^2$ and $0.005^2$ for reduced scattering and absorption, respectively. \\
\indent In contrast to the previous MAP estimates, by observing Figure \ref{fig:alternative_target}, we can see that the inverse-gamma and exponential hyperpriors have a substantial effect on the noise artefacts and sharpness of inclusions. Now, the inverse-gamma hypermodel produces sharp inclusions for both reduced scattering and absorption, while the exponential hyperprior greatly increases the contrast of the inclusions. On the other hand, the visibility of the background noise is also increased, which was expected due to the exponential hyperpriors feature to promote more common outliers. The standard gamma hyperprior produces sharp-edged

\noindent inclusions but, as a drawback, yields two false cross-talk inclusions for the reduced scattering. The cross-talk artefacts are also appearing in the other reconstructions, but with reduced contrast. \\
\indent In general, we observed that the exponential and inverse-gamma hyperprior worked well in cases with only positive (or negative) inclusions with a relatively large contrast. On the other hand, the estimates were more prone to produce cross-talk artefacts with the high contrast targets. \\

\subsection{Convergence of the nonlinear IAS}

\begin{figure}[t!]
\noindent\makebox[\textwidth]{\includegraphics[trim={1.8cm 0cm 0cm 0cm},clip,scale=0.38]{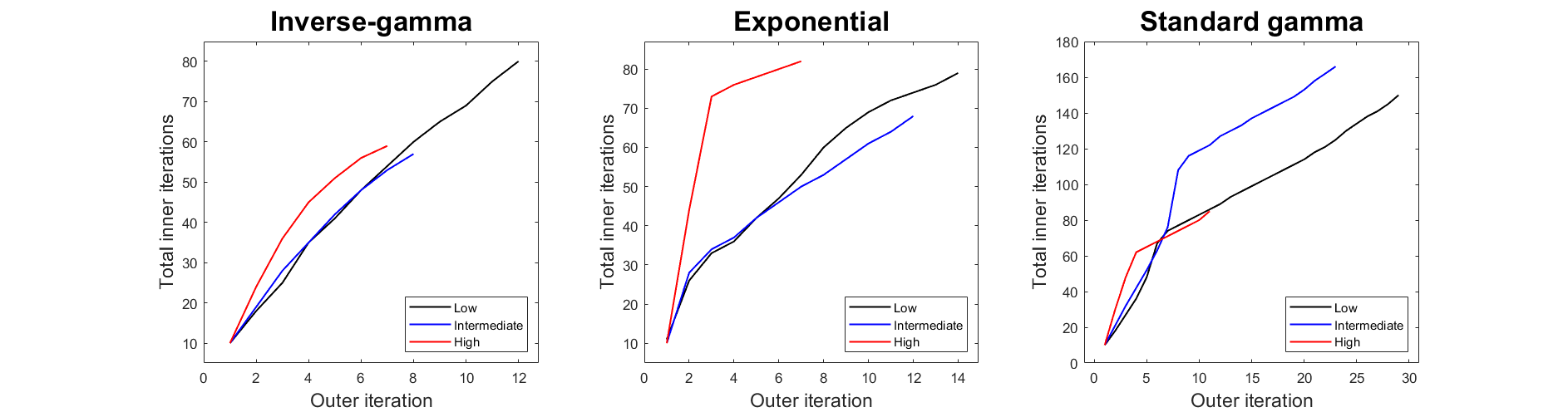}}
\caption{\revision{Computed inner (Gauss-Newton) iterations cumulatively as a function of the IAS iterations. The computed iterations were performed during the computation of the reconstruction with the uncorrelated Gaussian model as shown in Figure \ref{fig:fig_wn}.}}
\label{fig:inner_its}
\end{figure}

\begin{figure}[t!]
\noindent\makebox[\textwidth]{\includegraphics[trim={0 0.5cm 3.5cm 1cm},clip,scale=0.38]{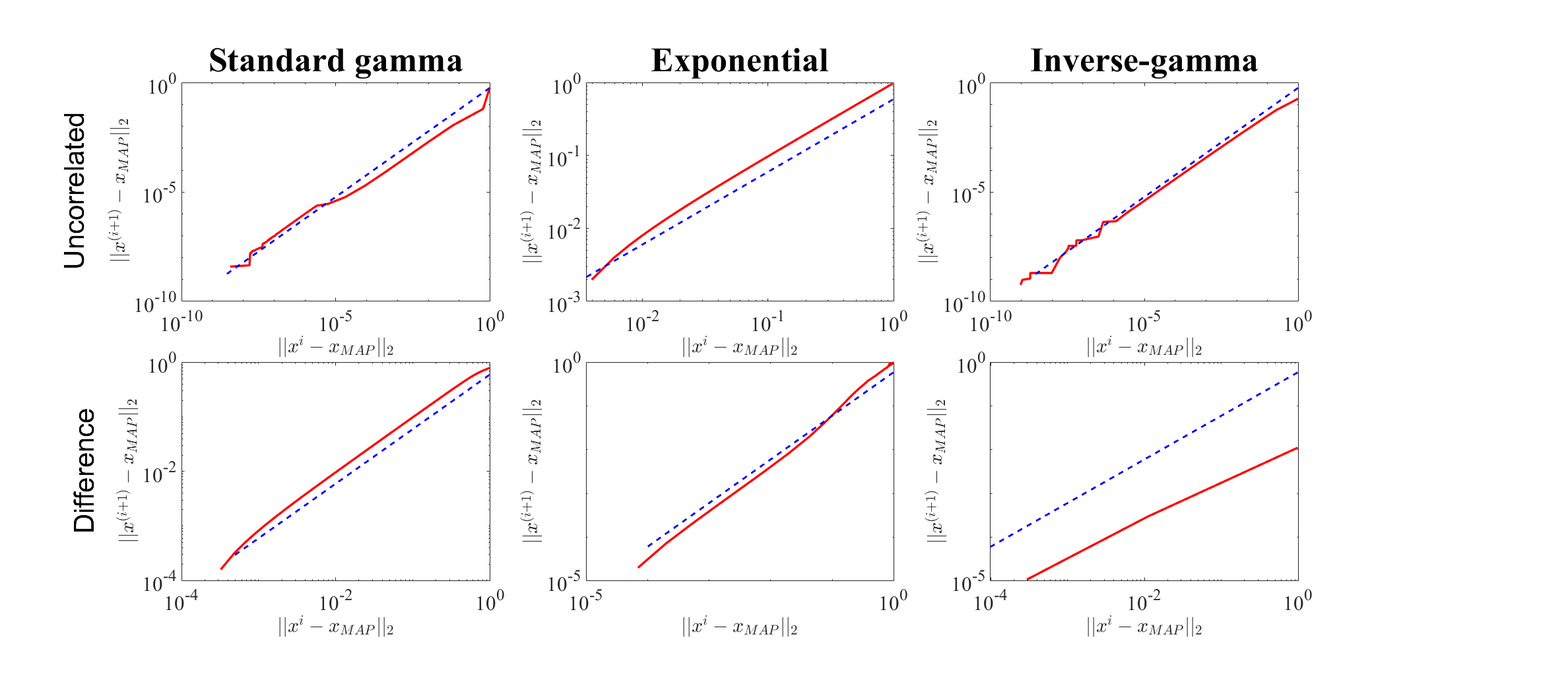}}
\caption{Convergence of the MAP estimates of the hierarchical models. The vertical axis shows the error $||x^i-x_{MAP}||_2$ of the $i$th iteration and the horizontal axis the error of the $(i+1)$th iteration. The MAP estimates $x_{\rm MAP}$ were computed by accurately solving the corresponding optimization problems. The bottom row shows the convergence of the difference prior model and the top row convergence of the uncorrelated Gaussian prior. The left, middle and right columns show the convergence for standard gamma, exponential and inverse-gamma hyperpriors, respectively. The red lines show the convergence of the used nonlinear IAS algorithm and the dashed blue line shows a linear convergence with convergence constant of $\mu=0.6$.} 
\label{fig:convergence}
\end{figure}

\begin{figure}[tbp!]
\noindent\makebox[\textwidth]{\includegraphics[trim={0 0cm 0 0},clip,scale=0.66]{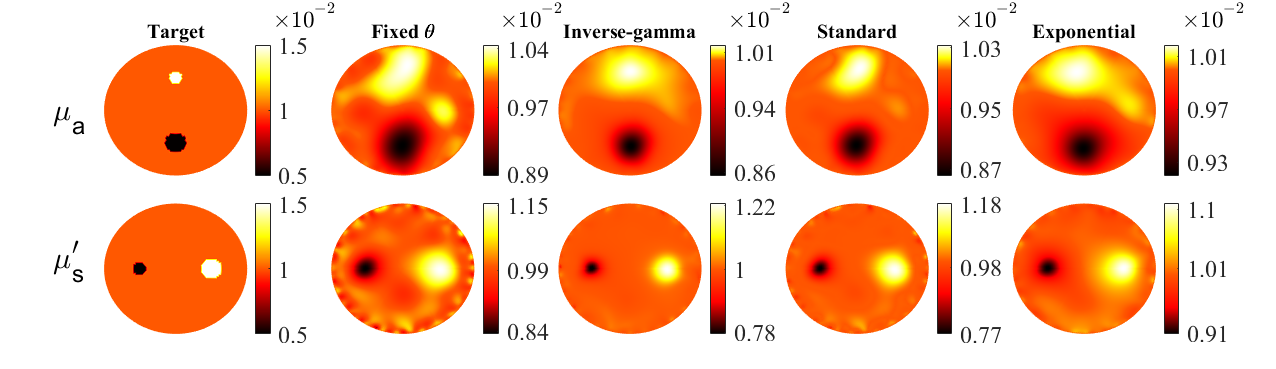}}
\caption{The absorption and reduced scattering estimates after three IAS iterations, when uncorrelated Gaussian prior was used with the intermediate hyperparameter values as in Figure \ref{fig:fig_wn}. }
\label{fig:early_grey}
\end{figure}

Besides the effect of the hierarchical models on the MAP estimates, we also wanted to investigate the empirical convergence of the nonlinear IAS \revision{and how the inner iterations behave. We restrict the convergence investigation to the case of uncorrelated Gaussian prior. \\
\indent A single Gauss-Newton iteration took 20-30 seconds for the models with the uncorrelated Gaussian prior and 40-60 seconds with the difference prior. Updating part b) of (\ref{eq:IAS_wh}) and \eqref{eq:ias_smoothness}, i.e., the variances during the IAS iterations took less than 0.1 seconds, which was neglectable compared to the total running time of a single IAS iteration. \\
\indent Figure \ref{fig:inner_its} shows the number of inner iterations over nonlinear IAS iterations cumulatively when using uncorrelated Gaussian prior. From Figure \ref{fig:inner_its}, we see that the first IAS iterations generally took more inner iterations, while this was reduced to just a few in later IAS iterations. By comparing the \textit{low}, \textit{intermediate} and \textit{high} curves, we can notice that the number of required IAS iterations tends to increase when the penalty term induced by the hyperprior is weighted more. On the other hand, when fewer IAS iterations were needed, the IAS iterations required more inner iterations (sharper rises in Figure  \ref{fig:inner_its}). The difference prior models were observed to need more IAS iterations than those with the uncorrelated Gaussian prior, while the required inner iterations had less deviation between IAS iterations.\\
\indent The inner iterations of the difference prior model followed essentially the same trends as with the uncorrelated Gaussian prior. Nevertheless, the average number of inner iterations per IAS iteration was much higher than with the uncorrelated Gaussian prior. This is to be expected since we also needed to approximately solve the compatibility condition during each step.}\\
\indent The convergence was observed for each MAP estimate \revision{of uncorrelated Gaussian model} with the intermediate hyperparameter values. That is, Figure \ref{fig:convergence} shows the convergence of the IAS, related to the MAP estimates shown in the middle columns of Figures \ref{fig:fig_wn} and \ref{fig:fig_smooth}. We also plotted a
linear convergence (blue dashed line) for reference in Figure \ref{fig:convergence}. By looking at the convergence slopes, the convergence of the nonlinear IAS can be observed to be comparable with the linear convergence. The linear convergence has a convergence constant, i.e., $\norm{x^{i+1}-x_{\rm MAP}}_2\leq\mu\norm{x^{i}-x_{\rm MAP}}_2$ with $\mu$ of 0.6. Thereby, the convergence of the nonlinear IAS was observed to be relatively fast in the linear sense. These empirical results are also in line with the previously proven \cite{calvetti19} at-least linear convergence for the uncorrelated Gaussian noise with standard gamma hyperpriors when the forward model is linear.\\
\indent Two important factors affecting the time-wise convergence are the required accuracy of the inner iterations (of Gauss-Newton) and the stopping criteria of the outer iterations. In this work, we used strict stopping conditions for the inner iterations to achieve an accurate minimizer of part a). On the other hand, a less accurate optimization of part a) could produce a sufficient approximation of the minimizer, allowing faster computation of the IAS iterations. For instance, with the linear forward model in MEG \cite{calvetti2015_krylov_MEG}, the Krylov-subspace approximation of the least squares solution has been shown to provide sufficient approximation leading to significant speedup. \\
\indent During the experimented IAS runs, it was observed that the first few iterations already encapsulate the main effect of the hierarchical models, especially with the uncorrelated Gaussian prior. After the early iterations, the IAS iterations only seemed to increase the contrast of the discontinuity points. Using early stopping criteria can drastically reduce the run time of the IAS algorithm, yet one needs to develop a systematic way to determine the criteria. To demonstrate the possibility of using early stopping, Figure \ref{fig:early_grey} shows the estimates of the optical parameters just after three iterations when the intermediate hyperparameter values were used as in Figure \ref{fig:fig_wn}. Finding a systematic way to set the early stopping criteria was beyond the scope of this work but could provide substantial speedup. For the nonlinear forward model, this would likely require further assumptions on the model.

\section{Conclusions}\label{sec:disc_conc}
In this work, hierarchical priors were formulated in the Bayesian framework for the highly ill-posed and nonlinear problem of diffuse optical tomography. The studied hyperpriors included exponential, standard, and inverse-gamma hyperpriors that were used with the uncorrelated Gaussian and difference prior. The simulated DOT problem demonstrated the hyperpriors having a considerable effect on the reconstructions of the piecewise linear targets. While the plausibility of the reconstructed characteristics depended on the chosen prior, hyperprior, and hyperparameters.  \\
\indent Results with the uncorrelated Gaussian prior showed the hierarchical models improving the localization of the inclusions by diminishing most of the background noise. The standard gamma hyperprior was observed to perform most robustly, while the exponential and inverse-gamma hyperpriors could not localize absorption inclusion well. When the difference prior was utilized, the hierarchical models were shown to enhance the sharpness of the edges. Again the standard gamma hyperprior was observed to work most robustly, with the broadest set of hyperparameters. The exponential and inverse-gamma hyperpriors could only enhance the edges of the simpler targets. Notably, the hierarchical priors showed excellent performance in supressing cross-talk artefacts.\\
\indent For all studied hierarchical models, the convergence of the nonlinear IAS was empirically investigated and was observed to be linear. These observations hold only when the inner (Gauss-Newton) iterations have strict stopping criteria, which can lead to slow computation times of a single IAS iteration. A computational speedup could be achieved by solving the inner iterations less accurately.\\
\indent In this work, we utilized the cumulative distribution approach to obtain suitable hyperparameter values for the hierarchical models with nonlinear forward operators. This method worked sufficiently well, but finding close to optimal hyperparameters still needed manual adjustment. Additionally, the cumulative distribution function approach excludes the possibility of considering spatial sensitivity differences, that is, a sensitivity-weighting. To justify more rigorous use of the hierarchical models with a nonlinear forward model, future work needs to focus on developing hyperparameter selection rules.

\section*{Acknowledgements}

This project has received funding from the European Research Council (ERC) under the European Union's Horizon 2020 research and innovation programme (grant agreement No 101001417 -- QUANTOM). 
This work was supported by the Finnish Cultural Foundation, the Academy of Finland (Center of Excellence in Inverse Modeling and Imaging project 336799, 336796, the Flagship Program Photonics Research and Innovation grant 320166, and the Academy Research Fellow project 338408, 346574)


\bibliographystyle{plain}

\bibliography{sources.bib}

\begin{thebibliography}{10}

\bibitem{aguerrebere2017bayesian}
C.~Aguerrebere, A.~Almansa, J.~Delon, Y.~Gousseau, and P.~Mus{\'e}.
\newblock A {B}ayesian hyperprior approach for joint image denoising and
  interpolation, with an application to hdr imaging.
\newblock {\em IEEE Transactions on Computational Imaging}, 3(4):633--646,
  2017.

\bibitem{arridge1999optical}
S.~R. Arridge.
\newblock Optical tomography in medical imaging.
\newblock {\em Inverse Problems}, 15(2):R41, 1999.

\bibitem{arridge2009}
S.~R. Arridge and J.~Schotland.
\newblock Optical tomography: forward and inverse problems.
\newblock {\em Inverse Problems}, 25:123010, 2009.

\bibitem{arridge1993performance}
S.~R. Arridge, M.~Schweiger, M.~Hiraoka, and D.~Delpy.
\newblock Performance of an iterative reconstruction algorithm for
  near-infrared absorption and scatter imaging.
\newblock In {\em Photon Migration and Imaging in Random Media and Tissues},
  volume 1888, pages 360--371. International Society for Optics and Photonics,
  1993.

\bibitem{bardsley2015randomize}
J.~M. Bardsley, A.~Seppanen, A.~Solonen, H.~Haario, and J.~P. Kaipio.
\newblock Randomize-then-optimize for sampling and uncertainty quantification
  in electrical impedance tomography.
\newblock {\em SIAM/ASA Journal on Uncertainty Quantification},
  3(1):1136--1158, 2015.

\bibitem{calvetti_MEG}
D.~Calvetti, H.~Hakula, S.~Pursiainen, and E.~Somersalo.
\newblock Conditionally {G}aussian hypermodels for cerebral source
  localization.
\newblock {\em SIAM Journal on Imaging Sciences}, 2(3), 2009.

\bibitem{calvetti2015_krylov_MEG}
D.~Calvetti, A.~Pascarella, F.~Pitolli, E.~Somersalo, and B.~Vantaggi.
\newblock A hierarchical {K}rylov-{B}ayes iterative inverse solver for {MEG}
  with physiological preconditioning.
\newblock {\em Inverse Problems}, 31(12):125005, 2015.

\bibitem{calvetti2020sparse}
D.~Calvetti, M.~Pragliola, E.~Somersalo, and A.~Strang.
\newblock Sparse reconstructions from few noisy data: analysis of hierarchical
  {B}ayesian models with generalized gamma hyperpriors.
\newblock {\em Inverse Problems}, 36(2):025010, 2020.

\bibitem{bootstrap}
D.~Calvetti, F.~Sgallari, and E.~Somersalo.
\newblock Image inpainting with structural bootstrap priors.
\newblock {\em Image and Vision Computing}, 24(7):782--793, 2006.

\bibitem{calvetti_uniform}
D.~Calvetti and E.~Somersalo.
\newblock Local regularization and {B}ayesian hypermodels.
\newblock {\em Proc. SPIE 5910, Advanced Signal Processing Algorithms,
  Architectures, and Implementations XV}, 5910, 2005.

\bibitem{calvetti2007gaussian}
D.~Calvetti and E.~Somersalo.
\newblock A {G}aussian hypermodel to recover blocky objects.
\newblock {\em Inverse problems}, 23(2):733, 2007.

\bibitem{calvetti2008hypermodels}
D.~Calvetti and E.~Somersalo.
\newblock Hypermodels in the {B}ayesian imaging framework.
\newblock {\em Inverse Problems}, 24(3):034013, 2008.

\bibitem{calvetti_bayes_overview}
D.~Calvetti and E.~Somersalo.
\newblock Inverse problems: From regularization to {B}ayesian inference.
\newblock {\em WIREs Computatinal Statistics}, 10(3), 2018.

\bibitem{calvetti19}
D.~Calvetti, E.~Somersalo, and A.~Strang.
\newblock Hierachical {B}ayesian models and sparsity: $\ell^2$-magic.
\newblock {\em Inverse Problems}, 35(3):035003, 2019.

\bibitem{chen2012optical}
J.~Chen.
\newblock Optical tomography in small animals with time-resolved {M}onte
  {C}arlo methods.
\newblock {\em Dissertation Abstracts International}, 74(04), 2012.

\bibitem{diamond2006}
S.~G. Diamond, T.~J. Huppert, V.~Kolehmainen, M.~A. Franceschini, J.~P. Kaipio,
  S.~R. Arridge, and D.~A. Boas.
\newblock Dynamic physiological modeling for functional diffuse optical
  tomography.
\newblock {\em Neuro{I}mage}, 30:88--101, 2006.

\bibitem{durduran2015}
T.~Durduran, R.~Choe, W.~B. Baker, and A.~G. Yodh.
\newblock Diffuse optics for tissue monitoring and tomography.
\newblock {\em Reports on Progress in Physics}, 73:076701, 2015.

\bibitem{gibson2005recent}
A.~Gibson, J.~Hebden, and S.~R. Arridge.
\newblock Recent advances in diffuse optical imaging.
\newblock {\em Physics in Medicine \& Biology}, 50(4):R1, 2005.

\bibitem{grosenick2016}
D.~Grosenick, H.~Rinneberg, R~Cubeddu, and P.~Taroni.
\newblock Review of optical breast imaging and spectroscopy.
\newblock {\em Journal of Biomedical Optics}, 21(10):091311, 2016.

\bibitem{guven2005diffuse}
M.~Guven, B.~Yazici, X.~Intes, and B.~Chance.
\newblock Diffuse optical tomography with a priori anatomical information.
\newblock {\em Physics in Medicine \& Biology}, 50(12):2837, 2005.

\bibitem{guven2005hierarchical}
M.~Guven, B.~Yazici, X.~Intes, and B.~Chance.
\newblock Hierarchical {B}ayesian algorithm for diffuse optical tomography.
\newblock In {\em 34th Applied Imagery and Pattern Recognition Workshop
  (AIPR'05)}, pages 6--pp. IEEE, 2005.

\bibitem{hebden2002three}
J.~C. Hebden, A.~Gibson, R.~Md. Yusof, N.~Everdell, E.~M.~C. Hillman, D.~T.
  Delpy, S.~R. Arridge, T.~Austin, J.~H. Meek, and J.~S. Wyatt.
\newblock Three-dimensional optical tomography of the premature infant brain.
\newblock {\em Physics in Medicine \& Biology}, 47(23):4155, 2002.

\bibitem{hiltunen2009combined}
P.~Hiltunen, S.~Prince, and S.~R. Arridge.
\newblock A combined reconstruction-classification method for diffuse optical
  tomography.
\newblock {\em Physics in Medicine \& Biology}, 54(21):6457, 2009.

\bibitem{hoshi2016}
Y.~Hoshi and Y.~Yamada.
\newblock Overview of diffuse optical tomography and its clinical applications.
\newblock {\em Journal of Biomedical Optics}, 21(9):091312, 2016.

\bibitem{isaacson1986distinguishability}
David Isaacson.
\newblock Distinguishability of conductivities by electric current computed
  tomography.
\newblock {\em IEEE transactions on medical imaging}, 5(2):91--95, 1986.

\bibitem{Ishimaru}
A.~Ishimaru.
\newblock {\em Wave Propagation and Scattering in Random Media}.
\newblock Academic Press, 1978.

\bibitem{Kaipio}
J.~P. Kaipio and E.~Somersalo.
\newblock {\em Statistical and Computational Inverse Problems}.
\newblock Springer, Newyork, 2005.

\bibitem{sens_2}
F.~Lin, T.~Witzel, S.~P. Ahlfors, S.~M. Stufflebeam, J.~W. Belliveau, and M.~S.
  H$\Ddot{\text{a}}$m$\Ddot{\text{a}}$l$\Ddot{\text{a}}$inen.
\newblock Assessing and improving the spatial accuracy in {MEG} source
  localization by depth-weighted minimum-norm estimates.
\newblock {\em Neuro{I}mage}, 31(1):160--171, 2006.

\bibitem{liu2018image}
S.~Liu, J.~Jia, D.~Yimi Y.~Zhang, and Y.~Yang.
\newblock Image reconstruction in electrical impedance tomography based on
  structure-aware sparse {B}ayesian learning.
\newblock {\em IEEE transactions on medical imaging}, 37(9):2090--2102, 2018.

\bibitem{mattout2006meg}
J.~Mattout, C.~Phillips, W.~Penny, M.~Rugg, and K.~Friston.
\newblock {MEG} source localization under multiple constraints: an extended
  {B}ayesian framework.
\newblock {\em Neuro{I}mage}, 30(3):753--767, 2006.

\bibitem{miyamoto2011phase}
A.~Miyamoto, K.~Watanabe, K.~Ikeda, and M.~Sato.
\newblock Phase diagrams of a variational {B}ayesian approach with {ARD} prior
  in {NIRS}-{DOT}.
\newblock In {\em The 2011 International Joint Conference on Neural Networks},
  pages 1230--1236. IEEE, 2011.

\bibitem{mozumder2021model}
M.~Mozumder, A.~Hauptmann, I.~Nissil{\"a}, S.~R. Arridge, and T.~Tarvainen.
\newblock A model-based iterative learning approach for diffuse optical
  tomography.
\newblock {\em IEEE Transactions on Medical Imaging}, 41(5):1289--1299, 2021.

\bibitem{mozumder2014compensation}
M.~Mozumder, T.~Tarvainen, J.P. Kaipio, S.~R. Arridge, and V.~Kolehmainen.
\newblock Compensation of modeling errors due to unknown domain boundary in
  diffuse optical tomography.
\newblock {\em JOSA A}, 31(8):1847--1855, 2014.

\bibitem{nummenmaa2007hierarchical}
A.~Nummenmaa, T.~Auranen, M.~H{\"a}m{\"a}l{\"a}inen,
  I.~J{\"a}{\"a}skel{\"a}inen, J.~Lampinen, M.~Sams, and A.~Vehtari.
\newblock Hierarchical {B}ayesian estimates of distributed {MEG} sources:
  theoretical aspects and comparison of variational and {MCMC} methods.
\newblock {\em Neuro{I}mage}, 35(2):669--685, 2007.

\bibitem{paulsen1996enhanced}
K.~Paulsen and H.~Jiang.
\newblock Enhanced frequency-domain optical image reconstruction in tissues
  through total-variation minimization.
\newblock {\em Applied optics}, 35(19):3447--3458, 1996.

\bibitem{Perona}
P.~Perona and J.~Malik.
\newblock Scale-space and edge detection using anisotropic diffusion.
\newblock {\em IEEE Transactions on Pattern Analysis and Machine Intelligence},
  12(7):629--639, 1990.

\bibitem{roininen2016hyperpriors}
L.~Roininen, M.~Girolami, S.~Lasanen, and M.~Markkanen.
\newblock Hyperpriors for mat\'{e}rn fields with applications in {B}ayesian
  inversion.
\newblock {\em arXiv preprint arXiv:1612.02989}, 2016.

\bibitem{TOAST}
M.~Schweiger and S.~R. Arridge.
\newblock The {T}oast++ software suite for forward and inverse modeling in
  optical tomography.
\newblock {\em Journal of Biomedical Optics}, 19(4):040801, 2014.

\bibitem{line_search}
M.~Schweiger, S.~R. Arridge, and I.~Nissil\"a.
\newblock Gauss–{N}ewton method for image reconstruction in diffuse optical
  tomography.
\newblock {\em Physics in Medicine \& Biology}, 50(10):2365, 2005.

\bibitem{shaw2014}
C.~B. Shaw and P.~K. Yalavarthy.
\newblock Performance evaluation of typical approximation algorithms for
  nonconvex {$\ell_p$}-minimization in diffuse optical tomography.
\newblock {\em Journal of the Optical Society of America A}, 31(4):852--862,
  2014.

\bibitem{shimokawa2012hierarchical}
T.~Shimokawa, T.~Kosaka, O.~Yamashita, N.~Hiroe, T.~Amita, Y.~Inoue, and
  M.~Sato.
\newblock Hierarchical {B}ayesian estimation improves depth accuracy and
  spatial resolution of diffuse optical tomography.
\newblock {\em Optics express}, 20(18):20427--20446, 2012.

\bibitem{Takeaki2013}
T.~Shimokawa, T.~Kosaka, O.~Yamashita, N.~Hiroe, T.~Amita, Y.~Inoue, and
  M.~Sato.
\newblock Extended hierarchical {B}ayesian diffuse optical tomography for
  removing scalp artifact.
\newblock {\em Biomedical Optics Express}, 4(11):2411--2432, 2013.

\bibitem{L1_methods}
D.~Vidaurre, C.~Bielza, and P.~Larra\ {n}aga.
\newblock A survey of {L}1 regression.
\newblock {\em International Statistical Review}, 81(3):361--387, 2013.

\bibitem{wheelock2019high}
M.~Wheelock, J.~Culver, and A.~Eggebrecht.
\newblock High-density diffuse optical tomography for imaging human brain
  function.
\newblock {\em Review of Scientific Instruments}, 90(5):051101, 2019.

\bibitem{yamashita2016multi}
O.~Yamashita, T.~Shimokawa, R.~Aisu, T.~Amita, Y.~Inoue, and M.~Sato.
\newblock Multi-subject and multi-task experimental validation of the
  hierarchical {B}ayesian diffuse optical tomography algorithm.
\newblock {\em Neuro{I}mage}, 135:287--299, 2016.

\bibitem{yoo2019deep}
Jaejun Yoo, Sohail Sabir, Duchang Heo, Kee~Hyun Kim, Abdul Wahab, Yoonseok
  Choi, Seul-I Lee, Eun~Young Chae, Hak~Hee Kim, Young~Min Bae, et~al.
\newblock Deep learning diffuse optical tomography.
\newblock {\em IEEE transactions on medical imaging}, 39(4):877--887, 2019.

\bibitem{zhang2012map}
G.~Zhang, X.~Cao, B.~Zhang, F.~Liu, J.~Luo, and J.~Bai.
\newblock {MAP} estimation with structural priors for fluorescence molecular
  tomography.
\newblock {\em Physics in Medicine \& Biology}, 58(2):351, 2012.

\end{thebibliography}

\end{document}